\documentclass[12pt]{article}
\usepackage{amsmath}
\usepackage{graphicx,psfrag,epsf}
\usepackage{enumerate}
\usepackage{natbib}

\def\e{\hbox{E}}

\def\var{\hbox{Var}}

\def\bias{\hbox{Bias}}
\def\se{\hbox{SE}}

\def\max{\hbox{max}}
\def\IF{\hbox{IF}}

\newcommand{\F}{\mathcal{F}}
\newtheorem{definition}{Definition}
\newtheorem{proposition}{Proposition}

\begin{document}

\def\spacingset#1{\renewcommand{\baselinestretch}%
{#1}\small\normalsize} \spacingset{1}


{
  \title{\bf Quantile Versions of the Lorenz Curve}
  \author{Luke A. Prendergast  and  Robert G. Staudte \\
    Department of Mathematics and Statistics, La Trobe University}
  \maketitle
} 

\bigskip
\begin{abstract}
The classical Lorenz curve is often used to depict inequality in a population of incomes, and the associated Gini coefficient is relied upon to make comparisons between
different countries and other groups. The sample estimates of these
 moment-based concepts are sensitive to outliers and so we investigate the extent to which quantile-based definitions  can capture income inequality and
 lead to more robust procedures. Distribution-free estimates of the corresponding coefficients of inequality are obtained, as well as sample sizes required to estimate them to a given accuracy. Convexity, transference and robustness of the measures are examined and illustrated.
\end{abstract}

\noindent%
\vfill
\newpage

 \section{Introduction}\label{sec:introduction}
The Lorenz curve and the associated Gini coefficient are routinely employed
for comparisons of income inequality in various countries. There are also numerous applications of them in the biological, computing, health and social sciences. These concepts have nice mathematical properties, and thus
are the subject of numerous theoretical  studies; for a recent review see \cite{kleiber-2005}.
However, when it comes to statistical inference for the Lorenz curve and the Gini coefficient, thorny issues arise.
An excellent review of existing methods and new proposals for estimating the standard error of the Gini coefficient
are investigated by \cite{davidson-2008}. However, as this author notes, such methods will not work when the variance
of the income distribution is large or fails to exist, and of course this means that they are undermined when
there are outliers in the data.
\cite{cowell-1996} show that most of the inequality measures in the econometrics literature
are very sensitive to outliers and have unbounded influence functions.

There are methods available for  resolving these inferential obstacles. One is to choose
a parametric income model and then to find optimal bounded influence estimators
for the parameters; for example, \cite{VFronc-1994} do this for the gamma and Pareto models. And, \cite{VF-2000} shows how to robustly choose between parametric models and then find robust estimates of inequality indices based on a single data sample, even if it has been grouped or truncated.
In a series of papers \cite{cowell-2002,cowell-2003,cowell-2007} investigate
damaging effects of data contamination on  transfer properties of various inequality indices, as well as dealing with the effects of truncation of non-positive and/or large data values. They propose semi-parametric models for overcoming these issues.

We go one step further here, redefining the basic concept of the Lorenz curve in terms of quantiles instead of moments, and then determining what has been gained and lost
in terms of conceptual clarity, inference and resistance to contamination.
Examples of this approach are the standardized median in lieu of the standardized mean, and quantile measures of skewness and kurtosis, rather than the classical
moment-based measures, \cite{S-2013,S-2014,S-2015}.  Ratios of quantiles
based on one sample are often presented as  measures of inequality, and inferential procedures for them are in \cite{ps-2015a,ps-2015b}.

The role of quantiles in inequality measures is long-standing.
\cite{gast-1971} observed that the definition of the Lorenz curve could be extended to all distributions having a finite mean $\mu $ by expressing the cumulative income as an integral of the quantile function.  More recently \cite{gast-2012} showed that  the inequality coefficient of \cite{gini-1914} could be made much more sensitive to shifts in income inequality if the mean in its denominator were replaced by the median; this also has the advantage of protecting the coefficient from large outliers.  \cite{kampke-2010} compares the effects of means versus
medians on poverty indices. It is in this spirit that we begin  in Section~\ref{sec:lorenzgini} by introducing three simple quantile versions of the Lorenz curve for distributions on the positive axis, and their associated coefficients of inequality.  Numerous examples demonstrate how these curves and coefficients agree or disagree with the moment-based classical version. In particular, the effects of an income transfer function on the inequality coefficients are illustrated for the Type~II Pareto model.

In Section~\ref{sec:finite} we study empirical versions of these inequality curves and their associated estimated coefficients. The latter estimates are found to have predictable
distribution-free standard errors, unlike the sensitive Gini coefficient. For an assumed
scale model, confidence intervals for the inequality coefficients are given.
It is not surprising that these quantile measures of inequality are resistant to outliers, and in Section~\ref{sec:robust} we show that they have bounded influence functions.

While the quantile versions of the Lorenz curve are not always convex, they are so for common distributions used to model incomes, as explained in Section~\ref{sec:convex}. A summary and further research problems are given in Section~\ref{sec:summary}.

\section{Quantile analogues of the Lorenz curve}\label{sec:lorenzgini}

\subsection{Definitions and basic properties}

Let $\F $ be the class of all cumulative distribution functions  $F$ with $F(0)=0.$  Such $F$ will be interpreted as \lq income\rq\ distributions and $p=F(x)$ as the
proportion of incomes less than or equal to $x.$
Define the {\em quantile function} associated with $F\in \F$ at each  $ p \in [0,1]$ by $Q(F;p)=F^{-1}(p)\equiv \inf \{x:\;F(x)\geq p\}$. If the support of $F$ is infinite; that is $F(x)<1$ for all $x>0$, this infimum does
not exist for $p=1$, and then we define $Q(F;1)=+\infty $.
When the meaning of $F$ is clear, we will sometimes write $x_p$ or $Q(p)$ for $Q(F;p)$.

  The mean income of those with proportion $p$  of  smallest incomes is $\mu =\mu _p(F)=\int _0^{x_p}x\,dF(x)/p$, and the mean income of the entire population  is defined by $\mu =\mu (F)=\lim _{p\to 1}\mu _p.$  Let $\F _0\subset \F$ be the set of $F$ for which $\mu (F)$ exists as a finite
 number. For each $F\in \F _0$  the  {\em Lorenz curve} of $F$ is defined by $L_0(F;p)\equiv p\,\mu _p/\mu ,$ for $0\leq p \leq 1.$  The lowest proportion of incomes $p$ have proportion $L_0 (p)$ of the total wealth.

 What we are proposing here is to replace $\mu _p$,
 the mean of the proportion $p$ of those with wealth less than $x_p,$ by its median $x_{p/2}=Q(F;p/2)$. In addition, we replace the mean $\mu $ of the entire population by  one of three quantile measures of its size:\quad $x_{1/2}$,
 $x_{1-p/2}$, or $(x_{p/2}+x_{1-p/2})/2.$ The robustness merits of this last divisor, a symmetric quantile average, are investigated by \cite{brown-1981}.

\begin{definition}\label{def1}
For $F\in \F$ and $p\in [0,1]$ let $x_{p}=Q(F;p).$
 The three  quantile-based functions whose graphs reveal income inequality are defined for each $p$ by:
 \begin{eqnarray}\label{qcurves} \nonumber
  L _1(F;p)&\equiv & p \,\frac {x_{p/2}}{x_{0.5}} \\
  L _2(F;p) &\equiv & p\, \frac {x_{p/2}}{x_{1-p/2}}\\ \nonumber
  L _3(F;p) &\equiv & 2p \,\frac{x_{p/2}}{ (x_{p/2}+x_{1-p/2})} = \frac{2p}{  1+p/L _2(F;p)}~.
\end{eqnarray}
As with $Q(p)=Q(F;p)$ we sometimes abbreviate $L_i(F;p)$ to $L_i(p)$.
\end{definition}

For each $p$ the first measure $L _1(p)$ compares the typical (median) wealth of the poorest proportion $p$ of incomes with the typical (median) wealth of the entire population. The second measure compares the bottom typical wealth with the top typical wealth; for example $L _2(0.2)$ corresponds to the popular \lq 20-20 rule\rq ,\ which compares the mean wealth  of the lowest 20\% of incomes with the largest 20\%. For each $p$ the third  $L _3(p)$ gives the typical wealth of the poorest $100\,p$\,\% incomes, relative to the mid-range wealth of the middle $100(1-p)$\,\% of incomes.
In all cases, extreme incomes are down-weighted because of multiplication by the factor $p$, as it is for the
Lorenz curve $L_0(p)=p\,\mu _p/\mu $.

\begin{figure}[t!]
\begin{center}
\includegraphics[scale=.5]{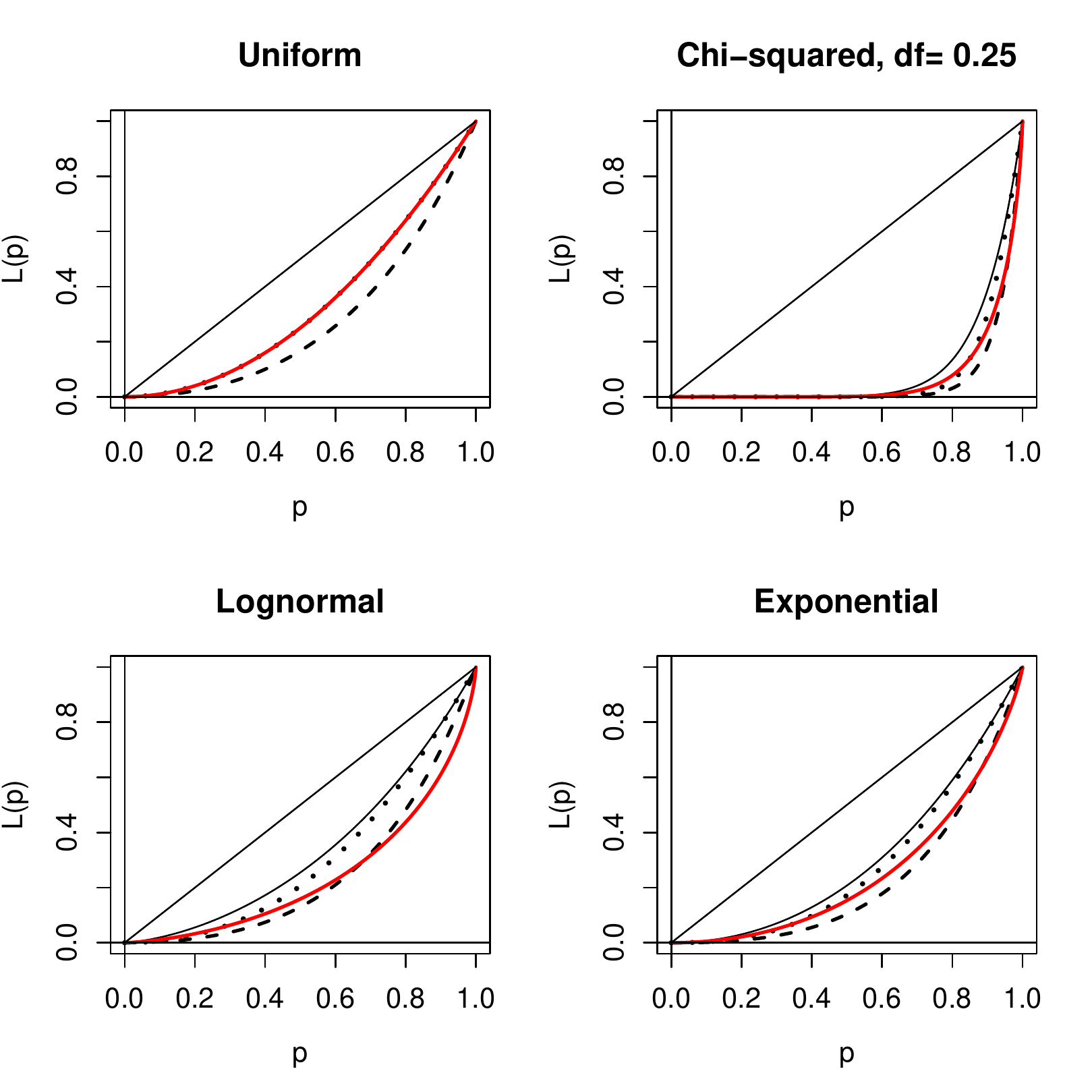}
\vspace{-.2cm}
\caption{\small \em Graphs of $L _1(p)$ (solid thin line),    $L _2(p)$ (dashed line),   $L _3(p)$ (dotted line), defined in (\ref{qcurves}) for various models. The thick solid  line is the Lorenz curve.  \label{fig1} }
\includegraphics[scale=.5]{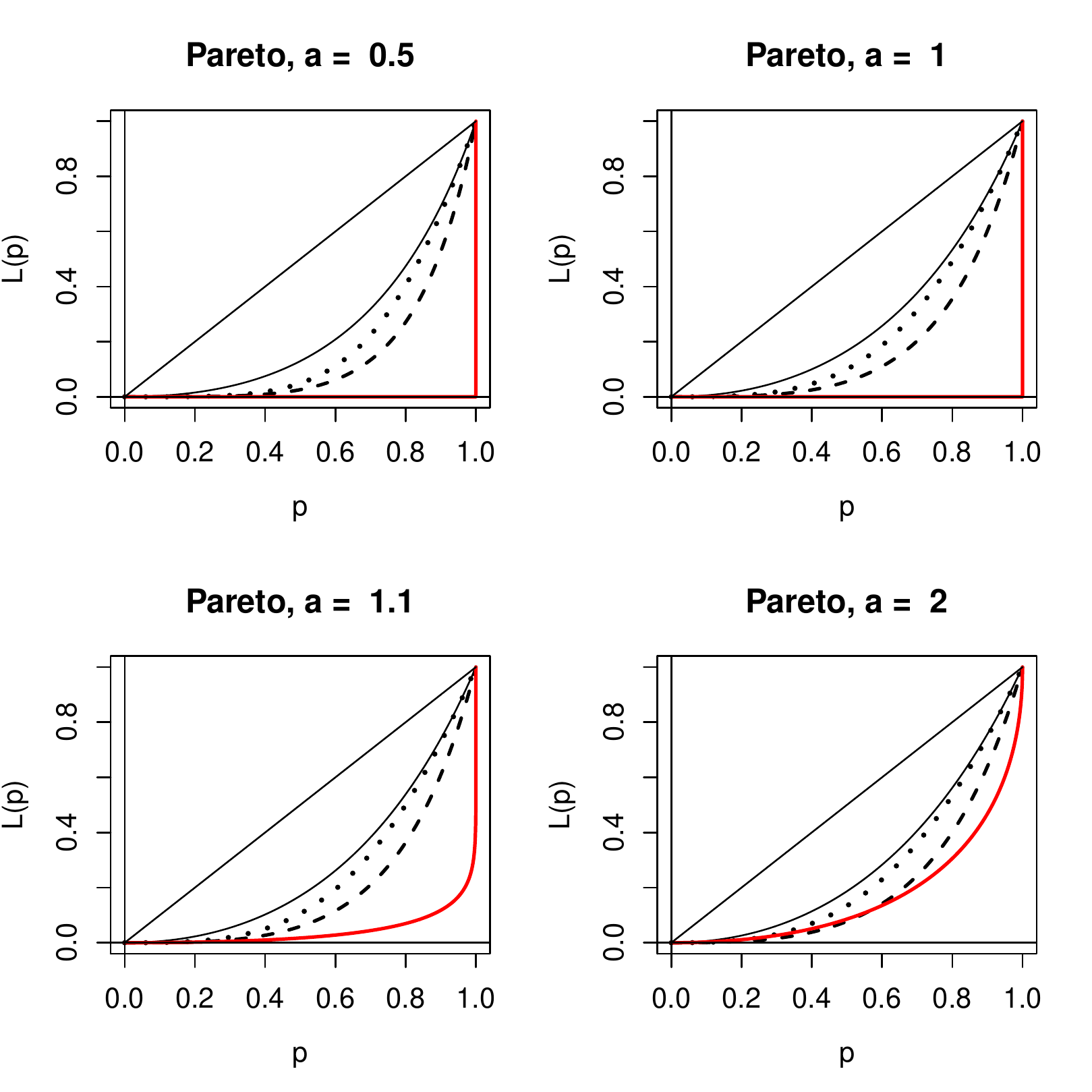}
\vspace{-.2cm}
\caption{\small Graphs of $L_i(p)$ for Type II Pareto$(a)$ Models, with the same notation as in Figure~\ref{fig1}.\label{fig2}}
\end{center}
\end{figure}

All of these quantile inequality curves $\{(p,L_i(p))\}$ are scale invariant and monotone increasing from $L_i(0)=0$  to $L_i(1)=1$, and all satisfy $L _i(p)\leq p$ for $0\leq p\leq 1$. Each $L_i(p)\equiv p$  when all incomes are equal. None are strictly speaking \lq Lorenz\rq\  curves, because they are not convex for all $F\in \F _0$, as  examples will show. Nevertheless, for most commonly assumed
income distributions $F$, they are convex, see Section~\ref{sec:convex}. Some examples of the quantile curves are depicted in Figures~\ref{fig1}-\ref{fig2}, which compares their graphs with the Lorenz curve. Note that $L_0(p)\equiv L_1(p)\equiv L_3(p)\equiv p^2$ for the uniform distribution.
And, $L_2(p)\approx p^3$ for the log-normal distribution.

\subsection{Coefficients of inequality}\label{QGini}

The relative measure of dispersion, or concentration ratio due to \cite{gini-1914} is defined for $F\in \F_0$ by $G_0=\e |X_1-X_2|/(2\mu ),$ where $X_1,X_2$ are independent and each distributed as $F$, and $\mu $ is the mean of $F$. It is known, see \cite{sen-1986} for example, to equal twice the area between the Lorenz curve and the diagonal line; it is an indicator, on the scale of 0 to 1, of \lq how far\rq\ the inequality graph is  from the diagonal line representing equal incomes; the further it is, the larger the Gini coefficient.

\clearpage
\newpage

\begin{definition}\label{def2}
For each of  the $L_i$ given in (\ref{qcurves}) define the respective {\em coefficients of inequality}
\begin{equation}\label{gini}
    G_i \equiv G _i(F)\equiv 2 \int _0^1\{p-L _i(F;p)\}dp\qquad \text { for all } F\in \F ~.
\end{equation}
\end{definition}
Specific numerical comparisons of the $G_i$s are given in Table~\ref{table1}.
It lists a variety of $F$ ranging from uniform to very long-tailed
 distributions and the associated values of Gini's index for the four $G_i$s. The  rankings of different $F$s by these four measures of inequality are seen to be very  similar and the Spearman rank correlation of $G_0$ with $G_i$  for $i=1,2$ and 3 are respectively 0.84, 0.88 and 0.88, for this list of $F$s. For more background material on distributions, see \cite{J-K-B-1994,J-K-B-1995}.
\begin{proposition}\label{P1} Given $F\in \F $ let $m=F^{-1}(0.5)$ denote its median.
Choose two incomes $Y_1,Y_2$ independently and randomly from those incomes less than the median, and let $V=\max \{Y_1,Y_2\}.$ It then follows that $G_1$ defined by (\ref{gini}) is given by $G_1 =\e [(m-V)/m],$ the relative average distance of
$V$ from the median. Next define $W=F^{-1}(1-F(V)),$ so if $V=x_r$ is the $r$th quantile of $F$, $W=x_{1-r}.$ It then follows that $G_2=\e [(W-V)/W]$ and $G_3=\e [(W-V)/(V+W)].$
\end{proposition}
{\bf Proof:\quad } {\em  Let $Y$ have the conditional distribution of $X$ given $X\leq m$; then its distribution function $F_Y(y)=2F(y),$ for $0\leq y\leq m$ and the distribution of $V$ is determined by $F_V(v)=F_Y^2(v)=4F^2(v)$, for $0\leq v\leq m.$ For each of  the three integrals in (\ref{gini}), make the change of variable $v=F^{-1}(p/2).$  The results are then immediate by observing that each  integral with respect to the measure $dF_V$ equals the corresponding claimed expected value.}

Proposition~\ref{P1} shows that $G_1\leq G_2$ and $G_3\leq G_2$ for all $F$.
It also allows for simple alternative interpretations of the three quantile inequality coefficients defined in (\ref{gini}) which can be compared with Gini's original definition as a relative measure of concentration.
Note that the Gini measure has been criticized for placing too much emphasis on
the central part of the distribution.  As Proposition~\ref{P1} shows, the quantile versions can also be criticized for the same reason, because the
main ingredient is the {\em maximum} of two randomly chosen incomes from the
lower half of the population. This maximum arises because in the definition (\ref{qcurves}) all the ratios are multiplied by $p$, which down-weights
ratios involving relatively small and large incomes.

\begin{table}[t!]
\begin{small}
\begin{center}
\caption{\label{table1} {\small \em Values of  $G _i$ to 3 decimal places for various $F$.} Also listed are the rankings of $F$ induced by the various $G _i$.}
\begin{tabular}{lrrrrrrrrrrrr}
 $\qquad F$ &&$G_0(F)$&$R_0$&& $G _1(F)$& $R_1$ && $G _2(F)$& $R_2$&& $G _3(F)$&$R_3$ \\
\hline
1. Uniform         && 0.333 &2          && 0.333 &4.5 && 0.455 &3  && 0.333 &3   \\
2. $\chi ^2_{0.5}$ && 0.762 &11         && 0.671 &13  && 0.792 &13 && 0.720 &13  \\
3. $\chi ^2_1$     && 0.636 &7          && 0.525 &11  && 0.673 &10 && 0.572 &10  \\
4. $\chi ^2_3$     && 0.423 &4          && 0.329 &3   && 0.483 &4  && 0.361 &4   \\
5. $\chi ^2_5$     && 0.339 &3          && 0.261 &2   && 0.406 &2  && 0.285 &2   \\
6. Lognormal       && 0.520 &6          && 0.333 &4.5 && 0.510 &5  && 0.388 &5   \\
7. Pareto(0.5)\footnotemark&& 1.000&13.5&& 0.515 &10  && 0.704 &11 && 0.610 &11 \\
8. Pareto(1)       && 1.000 &13.5       && 0.455 &9   && 0.636 &9  && 0.528 &9   \\
9. Pareto(1.5)     && 0.741 &9          && 0.434 &8   && 0.609 &8  && 0.497 &8   \\
10. Pareto(2)      && 0.667 &8          && 0.424 &7   && 0.595 &7  && 0.481 &7   \\
11. Weibull(0.25)  && 0.937 &10         && 0.731 &14  && 0.843 &14 && 0.787 &14   \\
12. Weibull(0.5)   && 0.750 &12         && 0.570 &12  && 0.720 &12 && 0.629 &12  \\
13. Weibull(1)     && 0.500 &5          && 0.393 &6   && 0.550 &6  && 0.432 &6   \\
14. Weibull(4)     && 0.159 &1          && 0.136 &1   && 0.222 &1  && 0.134 &1   \\[.2cm]
\hline
\end{tabular}
\end{center}
{\footnotesize $^{1.}$ The Lorenz curve and Gini coefficient
are not defined for distributions with $\mu = +\infty $, but if the definition were so extended,  $L_0(p)$ would be 0 for $0<p<1$ and the associated coefficient would be 1.}
\end{small}
\end{table}

\subsection{Tranference of income}\label{sec:transfer}

The effect of income transfers on inequality measures is of great interest to economists, see \cite{kleiber-2005} and \cite{fell-2012}. The basic idea \cite{dalton-1920} is that  if one transfers income from some members of the population having income above the mean to others having income below the mean, then the inequality measure should reflect this by decreasing. In keeping with our preference for quantiles
over moments, we suggest replacing the mean by the median in defining the transference principle for inequality measures.
\begin{definition}\label{def3}
Given $X\sim F\in \F $, and let $m\equiv x_{0.5}=F^{-1}(0.5)$ be the median.
We define a {\em median preserving transfer (of income) function} $Y=t(X)\sim F_Y$  as one
satisfying $m_y=F_Y^{-1}(0.5)=m$,
$F_Y(y)\leq F(y)$ for all $y<m$,  and $F_Y(y)\geq F(y)$ for all $y>m$.
\end{definition}
In words, a median preserving transfer function can only increase income that is less than the median, and only decrease income if it exceeds the median.
It follows that $y_p=Q(F_Y;p) \geq Q(F;p)=x_p$ for all $0<p<0.5$ and  $y_p=Q(F_Y;p) \leq Q(F;p)=x_p$ for all $0.5<p<1.$

The effect on the quantile inequality curves is then easily seen:\quad
$ L _1(F;p)=p \,x_{p/2}/x_{0.5}\leq p \,y_{p/2}/y_{0.5}=L _1(F_Y;p)$; that is, the transfer function
can only increase $L_1(p)$ at each $p$. This implies the associated coefficient of inequality (\ref{gini}) satisfies $G_1(F)\geq G_1(F_Y).$ We say that $L_1$ preserves the ordering induced by the transfer function. The reader may readily verify that for $i=2,3$
 the other quantile inequality curves satisfy $L _i(F;p)\leq L _i(F_Y;p)$ and hence $G_i(F)\geq G_i(F_Y).$

For any non-trivial transfer function we will have $G_i(F)>G_i(F_Y),$ a positive reduction in the coefficient of inequality.  Can we quantify this amount for any specific transfer functions? An example is given in Section~\ref{sec:extrans}.

\subsection{Example of transference}\label{sec:extrans}

Suppose one wants to increase all incomes less than a specific threshold $b$ (say the poverty line)
so that they equal $b$. That is; $t(x)=b$ for $0<x\leq b$. This requires an amount per person of
$d=b-(\int _0^bx\,dF(x))/F(b)$ to be found, say, by transference from those with incomes above the median or some higher thresh-hold $c$.  One possibility is to charge a levy of amount $d$ on those with income exceeding $c$,
leading to the following transfer function $Y=t(X)\sim F_Y$\ :
\begin{equation}\label{transferex}
    y=t(x)=\left\{
   \begin{array}{ll}
        b, & 0 \leq x<b ; \\
        x, & b\leq x < c; \\
        x-d, & c \leq x~.
             \end{array}
           \right.
\end{equation}
In the interest of fairness one could also charge a proportional  amount for those with income between $c$ and $c+d$ so that $Y=c$ for $c< x< c+d$,
but this unnecessarily complicates our presentation.

Now $F_Y(y)$ jumps from 0 to $F(b)$ at $b$, equals $F(y)$ for $b\leq y < c $, jumps at $c$ from
 $F(c)$ to $F(c+d)$ and equals $ F(y+d)$ for $ c \leq y$.
Therefore the quantile function $Q(F_Y;p)$ for the transferred income $Y$
is given by
\begin{equation}\label{Qtransferex}
    Q(F_Y;p)=\left\{
   \begin{array}{ll}
        b, & \qquad \quad \ 0\leq p<  F(b) ; \\
        F^{-1}(p), & \qquad F(b)\leq p <  F(c); \\
        c, &\qquad F(c) \leq p <  F(c+d); \\
        F^{-1}(p)-d, & \   F(c+d) \leq p ~.
             \end{array}
           \right.
\end{equation}

At this point it is convenient to introduce the {\em $p$th cumulative income} by $C(F;p)=\int _0 ^{x_p}y\,dF(y),$ where $x_p=Q(F;p)$. As \cite{cowell-2002} point out, this function is fundamental to analysis of Lorenz curves, and $C(1;F)=\mu $ and $L_0(F;p)=C(F;p)/C(1;F).$  We want to determine $C(F;p)$ for the Type II Pareto distribution having shape parameter $a>1$ and scale parameter $\sigma >0.$

Now $1-F_{a,\sigma }(x)= (1+x/\sigma )^{-a}$, which has mean $\mu = \sigma /(a-1)$ and $p$th quantile  $Q(F_{a,\sigma };\,p)=\sigma \{(1-p)^{-1/a}-1\}$.
Integrating by parts we obtain
\begin{equation}\label{cumincomePareto}
   C(F_{a,\sigma };\,p)= \int _0^{\sigma x_p}y\,dF_{a,\sigma }(y)=\frac {\sigma }{a-1}\left \{p-a(1-p)x_p\right \}~,
\end{equation}
where $x_p=Q(F_{a,1 };\,p).$ The mean income of the poorest proportion $p$ is $\mu _p =C(F_{a,\sigma };\,p)/p.$

 For the transfer problem with $F_{a,\sigma }(b)=p <0.5$, we have $b =\sigma \,x_{p }$, so (\ref{cumincomePareto}) implies
\[d =b -\mu _{p }=\frac {\mu }{p}\left \{(a-p) x_{p }-p )\right \}~.\]
This amount can  be obtained by a levy  $d $ on each income greater than $c =x_{1-p }.$

For the Pareto distribution with parameters $a=2$, $\sigma =100,000$ , the median income is 41,421.36 and the mean income is $\mu =100,000$.  For $p =0.2$, say, the quantities of interest are  the poverty line $b = 11,803.40$, the mean cumulative income $\mu _{0.2}= 5,572.80$ and $d =6,230.60.$    All those having income greater than the 0.8 quantile $123,606.30$ would need to pay an impost of $d =6,230.60$.

\begin{figure}[t!]
\begin{center}
\includegraphics[scale=.5]{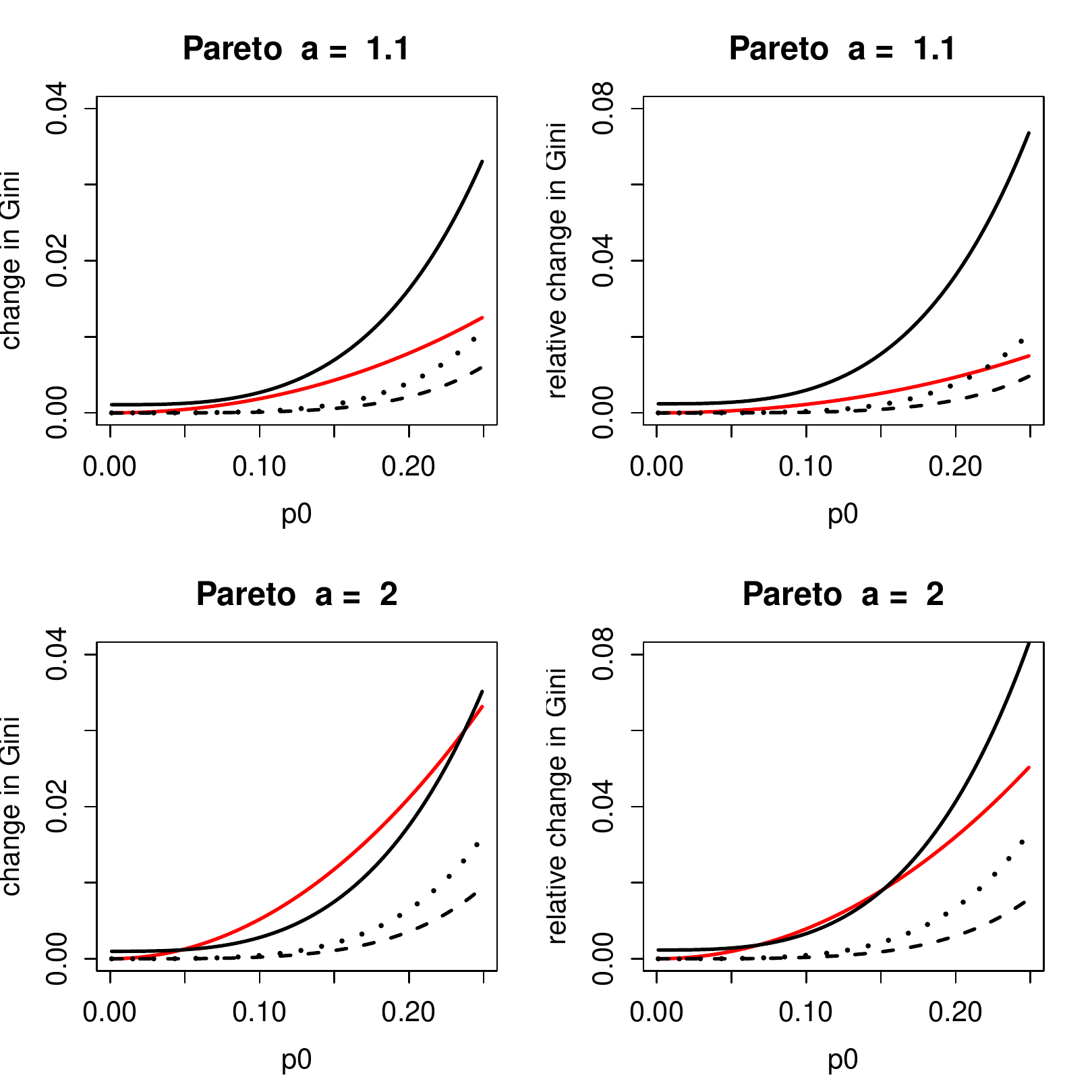}
\vspace{-.2cm}
\caption{\small \em The left hand plots show the graphs of the absolute change in the inequality coefficients $G_i(F;p)-G _i(F_Y;p)$ caused by the income
 transference (\ref{transferex}) for $i=0$, thick line; $i=1$, thin line; $i=2$, dashed  line; and $i=3$, dotted line. The right hand plots show the {\em relative} changes.\label{fig3}}
\end{center}
\end{figure}

The absolute and relative effects of such a transfer function are depicted in Figure~\ref{fig3} for two income distributions, Pareto with $a=1.1$ and $a=2$. For the first distribution, the change in the Gini coefficient $G_0$ is larger than for the $G_2$ and $G_3$ coefficients, but less than that for $G_1$; but the relative effect plot shows that the $G_1$ coefficient is most sensitive of the four, especially for $p_0$ near 0.25.
For the second distribution both $G_0$ and $G_1$ are roughly the same in terms of sensitivity to changes by transference and again preferable to $G_2$ and $G_3.$

Many other transfer functions and income distributions could be considered, but those are applications beyond the scope of this work.  It is important, of course, to identify real changes in coefficients of inequality after implementing a transfer of income. Estimation of $G_i$ is discussed in Section~\ref{sec:finite}. Another factor that we
have not included here are the costs of implementation of a transfer function.

\section{Estimation of inequality measures}\label{sec:finite}

\subsection{Empirical quantile inequality curves}\label{sec:finitelorenz}

Given data $x_1,\dots ,x_n$ with ordered values $x_{(1)}\leq x_{(2)} \leq \dots \leq x_{(n)}$ let $L_0(0)=0$ and $L_0(i/n)=\sum_{j\leq i} x_{(j)}/\sum_{j\leq n} x_{(j)}$ for $i=1,\dots ,n$. The empirical Lorenz curve is then defined as the graph of the piecewise linear connection of the points $(i/n,L_0(i/n))$, $i =0, 1,\dots ,n$.
The empirical distribution function defined for each $x$ by $F_n(x)=(\sum _{x_i\leq x}i)/n$ has inverse
$Q(F_n;p)=F_n^{-1}(p)=x_{([np]+1)}$ for $0\leq p <1$, and so empirical versions of the quantile curves (\ref{qcurves}) can be expressed in terms of the $n$ order statistics. Such curves are discontinuous, but there are several continuous quantile estimators available, including kernel density estimators \cite{shma-1990} and the linear combinations of two adjacent order statistics studied by \cite{hynd-1996}.  Many of the latter are implemented on the software package $R$ \cite{R}, and here we use the Type~8 version of the {\tt quantile} command recommended by  \cite{hynd-1996}. It linearly interpolates between the points $(p_{[k]},x_{(k)})$, where  $p_{[k]} = (k - 1/3) / (n + 1/3)$ and is a continuous function of $p$ in $(0,1).$
We also denote this estimator $\hat x_p=\hat Q(p)$.

\begin{definition}\label{def4}
All of the $L _i$ curves defined by (\ref{qcurves}) are functions of the quantile function $Q(F;p)$, so given the estimator $\hat x_p=\hat Q(p)$ one can by substitution obtain estimators of each
of the $L _i(p)$  for any  $p$ in $(0,1);$ we call these estimators $\hat L _i(p),$ for $i=1,2,$ and 3.
\end{definition}

\subsection{Empirical coefficients of inequality}\label{sec:finitegini}

With few exceptions, such as the uniform distribution, one cannot analytically compute the $G _i(F)$s, but using
modern software packages such as R \cite{R}, it is easy to get very good approximations to them for many $F$ of interest as follows.

\begin{definition}\label{def5}
Given a large integer $J$ define a grid in (0,1) with increments of size
$1/J$ by $p_j=(j-1/2)/J$, for $j=1,2,\dots ,J.$
 Then evaluate the quantile function $Q(p_j)$ for $p_j$ in the grid and find $G _i(J)\equiv  (2/J)\,\sum _j\{p_j -L _i(p_j)\}$ for each $i=1,2$ and 3.

Clearly one can make $G_i(J)$ as close to $G_i$ as desired by choosing $J $ sufficiently large. We will estimate $G_i(J)$, and hence $G_i$, as follows.
Let $\hat L_i(p_j)$ be the estimated inequality curve value at $p_j$, for
each $p_j$ in the grid. Then $\hat G _i(J)$ is defined by
\begin{equation}\label{ginihat}
    \hat G _i(J)\equiv   (2/J)\, \sum _j\{p_j-\hat L _i(p_j)\}~.
\end{equation}
\end{definition}

 In our computations, we used $J=1000$. Hereafter we write
$G_i$ for $G_i(J)$ and $\hat G_i$ for $\hat G_i(J ),$ but it
is understood that these are computed on a grid with increments $1/J.$

\begin{figure}[t!]
\begin{center}
\includegraphics[scale=.5]{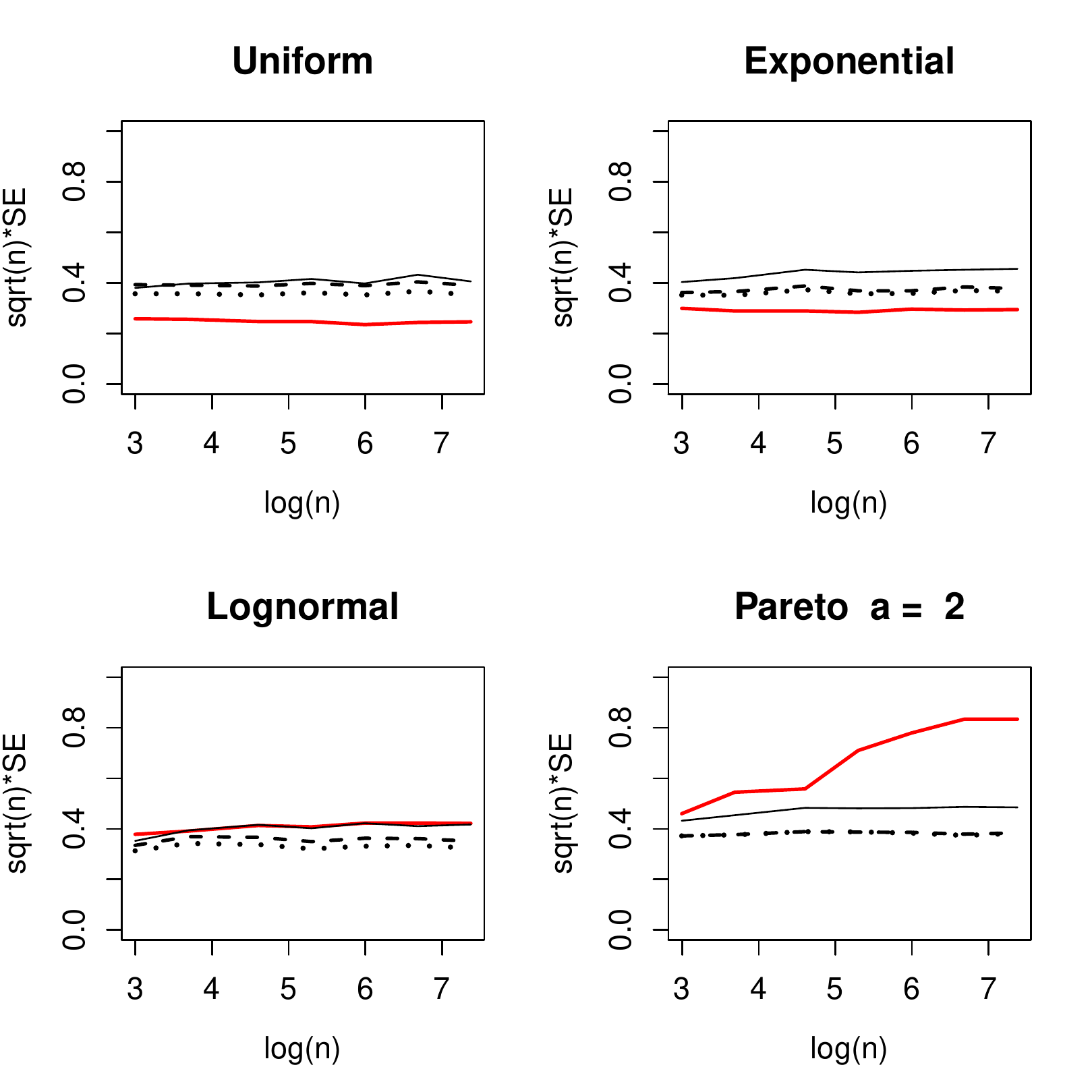}
\vspace{-.2cm}
\caption{\small \em  $\sqrt n\,\se [\hat G _i(F)]$ plotted as a function of $\ln (n)$ for the Lorenz curve $L_0$ (thick solid line) and
$L _i$-curves, (thin solid, dashed, and dotted lines,respectively). \label{fig4}}
\includegraphics[scale=.5]{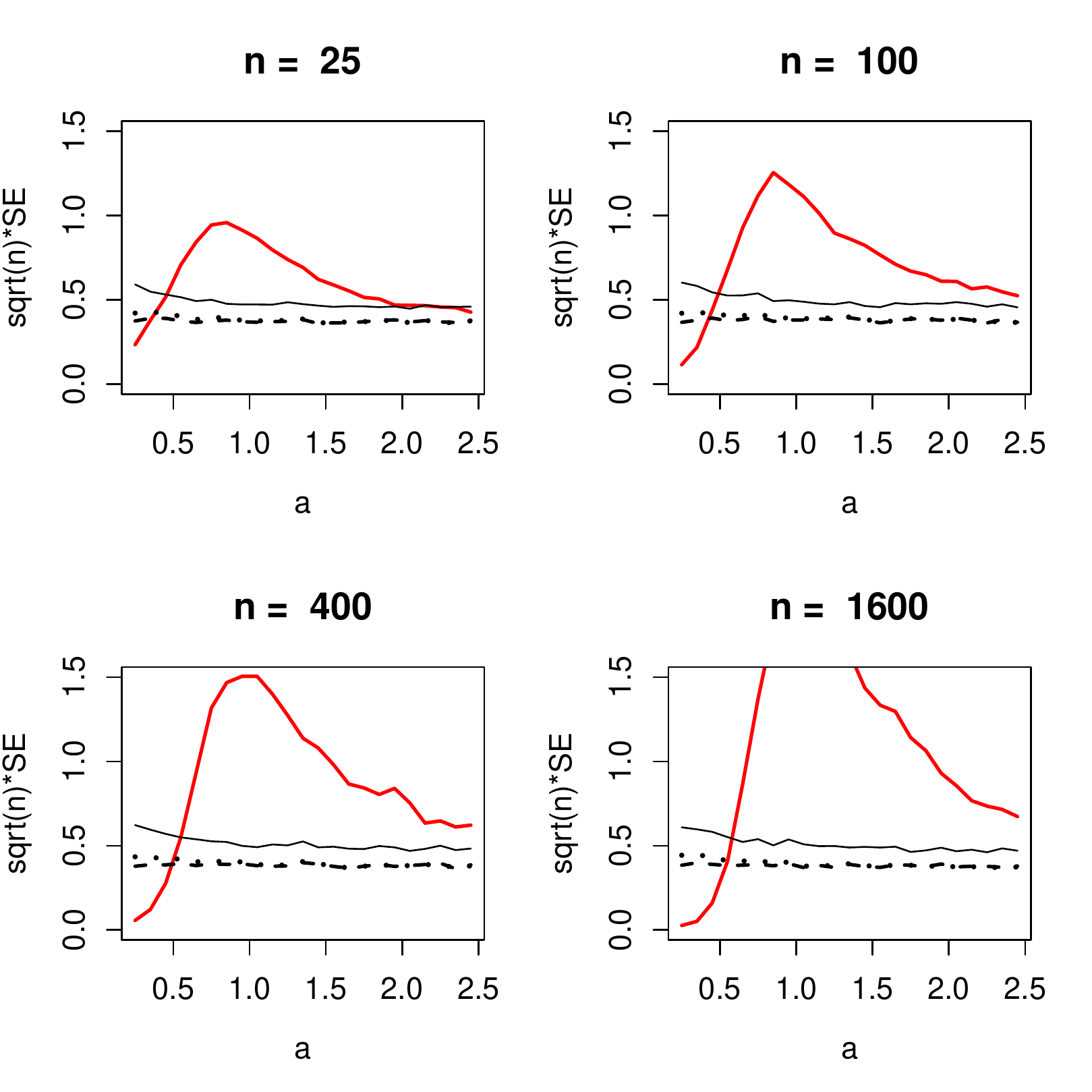}
\vspace{-.2cm}
\caption{\small \em $\sqrt n\,\se [\hat G _i(F_a)]$, for Pareto$(a)$ distributions, plotted as a function of $a$.\label{fig5}}
\end{center}
\end{figure}

\subsection{Simulation studies}\label{sec:simulations}

It will be seen that despite the fact that the values of the quantile coefficients of inequality $G _i(F)$ vary greatly over the range of $F$ in Table~\ref{table1}, the standard errors of estimation are fairly predictable. By \lq standard error\rq\  of $\hat G _i$, we mean the square root of the mean squared error.
Initial simulations suggested that $\bias [\hat G _i]=o(n^{-1/2})$ and $\var [\hat G _i]=O(1/n)$ so in Figure~\ref{fig4} we show some examples of
$\sqrt n\,\se [\hat G _i(F)]$, plotted as a function of $\ln (n),$ for $n$ ranging from
20 to 1600. These plots are based on 1000 replications at each of the selected values of $n$ for various $F$. In all four plots it is seen that the standard errors of $\hat G _2(F)\approx \hat G _3(F)\approx 1/(2\,\sqrt {n})$ while $\hat G _1(F)$ is a little larger.  This enables one to choose a sample size which guarantees a desired standard error for each of the three estimators.  Attempting to estimate Gini's coefficient of inequality by means of the Lorenz curve areas has no such simple sample size solution.

It is also interesting to plot $\sqrt n\,\se [\hat G _i(F_a)]$ versus $a$ as in Figure~\ref{fig5}, where $F_a$ denotes the Pareto distribution with shape parameter $a$
ranging from $0.25:2.5/0.1$.  Again all three standard errors of the estimated inequality coefficients derived from the $L _i$-curves are well behaved, but those for the Lorenz curve
are quite irregular. For $a\leq 1$ the Lorenz curve is not defined because $\e _a[X]=+\infty $ but if one defines the curve to be 0 in this case the corresponding measure of inequality
is 1 and this can be estimated. Even if one restricts attention to $1<a<2$, these plots show that for increasing $n$ the standard error is growing at a faster rate than the others, (because for $a<2$ the variance of $F$ is infinite).

The results in Table~\ref{table2} suggest that one can choose the minimum sample size required to obtain $\se [\hat G _1]\leq c$; it is $n_1(c)=(0.55/c)^2$. So for example, for standard error $c=0.01$, one needs $n\geq n_1\approx 3000$. {\em Note that this accuracy
is achieved for all $F$ in Table~\ref{table2}.}
Similarly for $G _2, G _3$ the required
sample size is a little smaller $n_2(c)=(0.43/c)^2=n_3(c).$

 \begin{table}[t!]
\begin{small}
\begin{center}
\begin{tabular}{lrrrrrrrrrrr}
 &\multicolumn{3}{c}{$\quad G _1$}& &\multicolumn{3}{c}{$G _2$}& &\multicolumn{3}{c}{$G _3$}\\
 $\qquad F$ &  25 & 100 &
  $+\infty $ & & 25 & 100  & $+\infty $ & & 25 & 100  & $+\infty $ \\
  \cline{2-4}\cline{6-8}\cline{10-12}
 1. Uniform         & 0.40 & 0.40 & 0.421 && 0.38 & 0.39 & 0.399 && 0.35 & 0.35 & 0.361\\
 2. $\chi ^2_{0.5}$ & 0.55 & 0.55 & 0.550 && 0.39 & 0.38 & 0.359 && 0.43 & 0.43 & 0.405\\
 3. $\chi ^2_1$     & 0.50 & 0.53 & 0.521 && 0.40 & 0.41 & 0.402 && 0.42 & 0.44 & 0.427\\
 4. $\chi ^2_3$     & 0.39 & 0.40 & 0.408 && 0.34 & 0.36 & 0.351 && 0.31 & 0.33 & 0.316\\
 5. $\chi ^2_5$     & 0.32 & 0.33 & 0.337 && 0.30 & 0.32 & 0.305 && 0.26 & 0.27 & 0.253\\
 6. Lognormal       & 0.39 & 0.40 & 0.417 && 0.34 & 0.35 & 0.351 && 0.32 & 0.32 & 0.322\\
 7. Pareto(0.5)     & 0.53 & 0.54 & 0.540 && 0.38 & 0.39 & 0.351 && 0.41 & 0.42 & 0.370\\
 8. Pareto(1)       & 0.49 & 0.50 & 0.507 && 0.37 & 0.38 & 0.371 && 0.38 & 0.39 & 0.376\\
 9. Pareto(1.5)     & 0.46 & 0.47 & 0.492 && 0.36 & 0.38 & 0.379 && 0.36 & 0.38 & 0.380\\
 10. Pareto(2)      & 0.45 & 0.46 & 0.485 && 0.37 & 0.38 & 0.381 && 0.37 & 0.38 & 0.379 \\
 11. Weibull(0.25)  & 0.55 & 0.53 & 0.540 && 0.35 & 0.34 & 0.330 && 0.40 & 0.39 & 0.384\\
 12. Weibull(0.5)   & 0.53 & 0.53 & 0.550 && 0.38 & 0.39 & 0.387 && 0.41 & 0.42 & 0.422\\
 13. Weibull(1)     & 0.44 & 0.45 & 0.461 && 0.37 & 0.38 & 0.382 && 0.36 & 0.37 & 0.370\\
 14. Weibull(4)     & 0.19 & 0.19 & 0.195 && 0.20 & 0.21 & 0.207 && 0.14 & 0.14 & 0.140
\end{tabular}
  \end{center}
  \end{small}
 \caption{\small {\bf Standard errors of $\hat G _i$}: \quad $\sqrt n\, \se [\hat G _i]$ for $n =25,100$ together with the respective asymptotic SEs $\sigma _i=\lim \sqrt n\, \se [\hat G _i]$, based on numerical evaluation of the integrals in (\ref{IFginimoments}). The finite sample  standard errors are based on 4000 samples of size $n$.\label{table2}}
 \end{table}

\subsection{Confidence intervals for the coefficients of inequality}\label{cigini}

Recall from (\ref{ginihat}) that for each $i=1,2,3$ and large fixed $J$ the estimated coefficient of inequality  is  $ \hat G _i=(2/J) \sum _j\{p_j-\hat L _i(p_j)\}.$  Now the estimate $\hat L _i(p_j)$, as a ratio of finite linear combinations of quantile estimates, is consistent for $ L _i(p_j)$, so $ \hat G _i$ is also consistent for $G _i$.  Further, \cite{ps-2015b} show that $n^{1/2}\{\hat L _i(p_j)-L _i(p_j)\}$ is asymptotically normal with mean 0 and variance depending on certain quantiles and quantile densities of the underlying $F$. \cite{beach-1983} find the limiting {\em joint} normal distribution of estimates of a finite number of Lorenz curve ordinates, based on a finite number of sample quantiles, assuming that $F\in \F _0\cap \F ',$
where $\F '$ is specified in Definition~\ref{def6}. In the same way, for $F\in \F '$,
the  limiting joint normal distribution of the estimated ordinates $\hat L_i(p_j)$, $j=1,\dots ,J$ can be established. We do not need an analytic expression for the covariance matrix, because we only require the asymptotic normality of the
estimated $G _i$, which being an average of the $p_j-\hat L _i(p_j)$, is immediate. Its asymptotic variance is available from the expected value of the squared influence function see (\ref{IFgini}) and (\ref{IFginimoments}).

Here we present the results of a modest simulation study of confidence intervals for $G _i$ of the form $\hat G_i\pm 1.96\sigma _i/\sqrt n\,$, with nominal coefficient 95\% , with results in Table~\ref{table3}.  For the lognormal distribution, the respective $\sigma _i$ found in Table~\ref{table2} are respectively 0.417, 0.351 and  0.322. For the Pareto with $a=2$ distribution, these values are 0.485, 0.381 and 0.379.

To obtain distribution-free confidence intervals for $G_i$, one needs consistent estimates for the asymptotic variance $\sigma _i^2$, a project beyond the scope
of this work.
\begin{table}[t!]
\begin{small}
\begin{center}
\begin{tabular}{lrrrrrrrrrrr}
 &\multicolumn{3}{c}{$\quad G _1$}& &\multicolumn{3}{c}{$G _2$}& &\multicolumn{3}{c}{$G _3$}\\
 $\qquad F$ &  25 & 100 & 400 & & 25 & 100  & 400 & & 25 & 100  & 400 \\
  \cline{2-4}\cline{6-8}\cline{10-12}
Lognormal  & 0.967  & 0.956 & 0.951 && 0.954 & 0.947 & 0.947 && 0.954 & 0.946 & 0.948\\
           & 0.327  & 0.164 & 0.082 && 0.275 & 0.138 & 0.069&& 0.252 & 0.126 & 0.063 \\
           & \\
Pareto(2)  &  0.966 & 0.960 & 0.956 && 0.955 & 0.952 & 0.951 && 0.954 & 0.951& 0.950  \\
   &  0.380 & 0.190 & 0.095 && 0.299 & 0.149 & 0.075&& 0.297 & 0.149& 0.074
\end{tabular}
  \end{center}
  \end{small}
 \caption{\small {\bf Confidence intervals for $ G _i$}: \quad \em Empirical coverage probabilities and widths based on 10,000 simulations of nominal  95\% confidence intervals.\label{table3} }
 \end{table}

\section{Robustness properties}\label{sec:robust}

In this section we show that the quantile inequality curves and their associated
coefficients of inequality have bounded influence functions, which guarantees
that a small amount of contamination can only have a limited effect on the asymptotic bias of estimators of these quantities. For background material on robustness concepts for functionals, see \cite{HRRS86}, although we attempt to
make the presentation self-contained. To this end, we must restrict $F\in \F$
to the following subclass of smooth distributions:

\begin{definition}\label{def6}
\begin{equation*}
   \F '=\{F\in \F :\ f=F' \text { exists and is strictly positive.} \}
\end{equation*}
For $F\in \F '$ with inverse $x_p=Q(p)=F^{-1}(p)$, we define the
 {\em quantile density} by
\begin{equation}\label{qden}
 q(p)=\frac {\partial \;Q(F;p)}{\partial \;p}\,=\, \frac {1}{F'(Q(F;p))}\,=\,\frac {1}{f(x_p)}~.
\end{equation}
The quantile density terminology is due to \cite{par-1979}, and its importance
was earlier recognized by \cite{tukey-1965} who called it the \lq sparsity index\rq .\
\end{definition}

In order to find the influence function of the $L _i$-curves at any specific $p$ in $(0,1)$ we also require  the mixture distribution which places positive  probability $\epsilon $ the point $z$ (the contamination point) and $1-\epsilon $ on the income distribution $F$. Formally, it is defined for each $x$ by $F_\epsilon ^{(z)}(x)\equiv (1-\epsilon)F(x)+\epsilon I[x\geq z]$, where $I[\cdot ]$ is the indicator function.  The {\em influence function} for any functional $T$ is then defined for each $z$ as  the $\IF (z;T,F)\equiv \lim _{\epsilon \downarrow 0}\{T( F_\epsilon ^{(z)})-T(F)\}/\epsilon =\lim _{\epsilon \downarrow 0}\frac {\partial }{\partial \epsilon }T( F_\epsilon ^{(z)})$.
 The influence function of the $p$th quantile functional $T(F)=Q(F;p)$, where $F\in \F '$ of Definition~\ref{def6}, is well-known to be \cite[p.59]{S-S-1990}
\begin{equation}\label{IFquantile}
    \IF (z;\,Q(\,\cdot\, ;p),F)=\left\{
                                  \begin{array}{ll}
                                 (p-1)\,q(p) , & \hbox{$z< x_p$\,;} \\
                                   0 , &  \hbox{$z= x_p$\,;} \\
                                  p\,q(p)  , & \hbox{$z> x_p$\,.}  \end{array}
                                \right.
\end{equation}
where $x_p=F^{-1}(p)$ and $q(p)$ is given by (\ref{qden}).

 One can show that
 $\e _F [\IF (Z;\,Q(\,\cdot\, ;p),F),F)]=0$ and $\var _F [\IF (Z;\,Q(\,\cdot\, ,p),F),F)]=\e _F [\IF ^2(Z;\,Q(\,\cdot\, ,p),F),F)]=p(1-p)\,q^2(p)$. One reason for calculating this variance is that it arises in the asymptotic variance of the functional applied to the empirical distribution $F_n$, namely $Q(F_n ;p)$. That is, $n^{1/2}\; [Q(F_n ;p)-Q(F ;p)]\to N(0,p(1-p)\,q^2(p))$ in distribution; and sometimes a simple expression for the asymptotic variance is not otherwise available.

\begin{figure}[t!]
\begin{center}
\includegraphics[scale=.5]{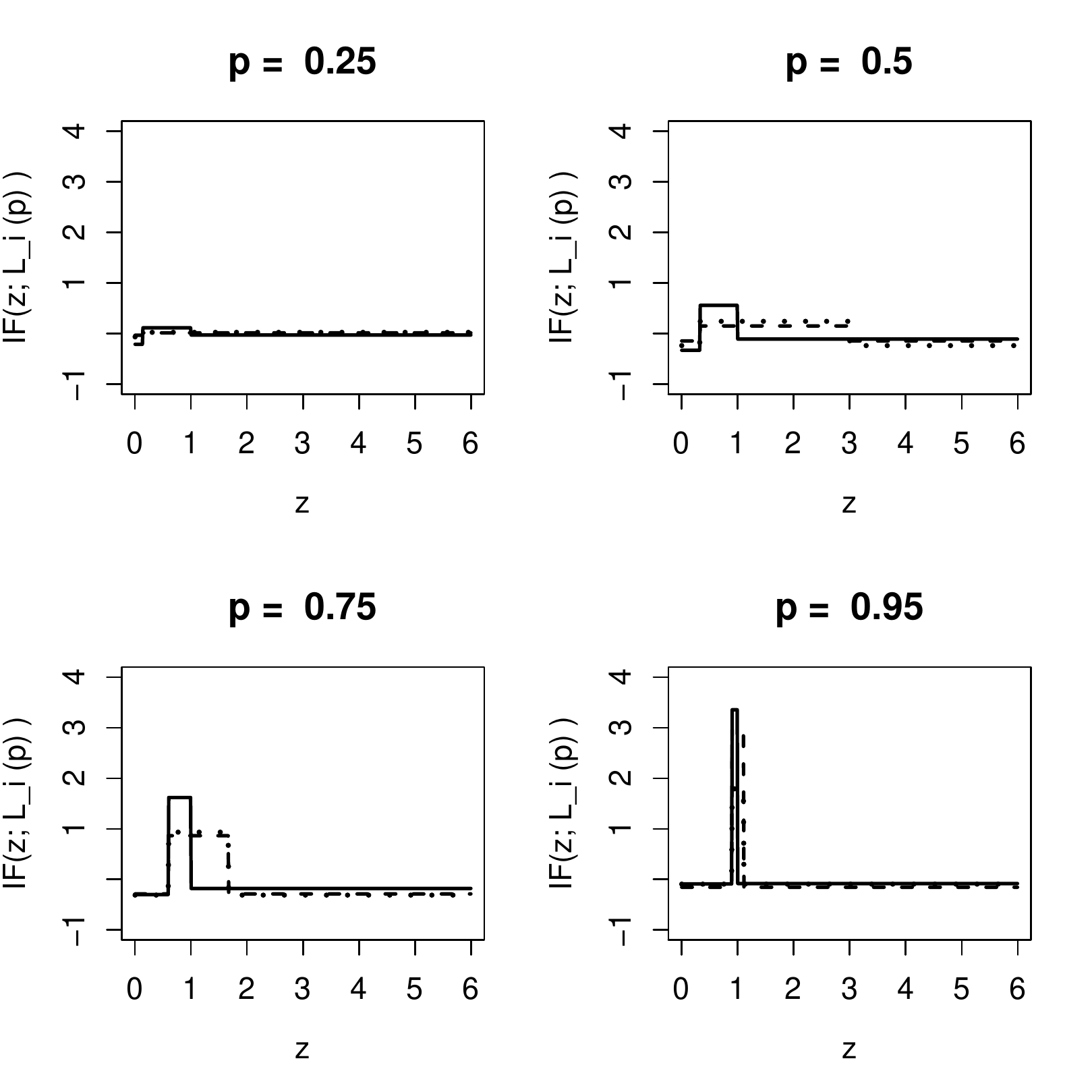}
\vspace{-.4cm}
\caption{\small \em For various choices of $p$, $\IF (z;\,L _i(p),F_1)$ is plotted as a function of $z$; the solid, dashed and dotted lines correspond, respectively, to $i=1,2$ and 3. \label{fig6}}
\includegraphics[scale=.5]{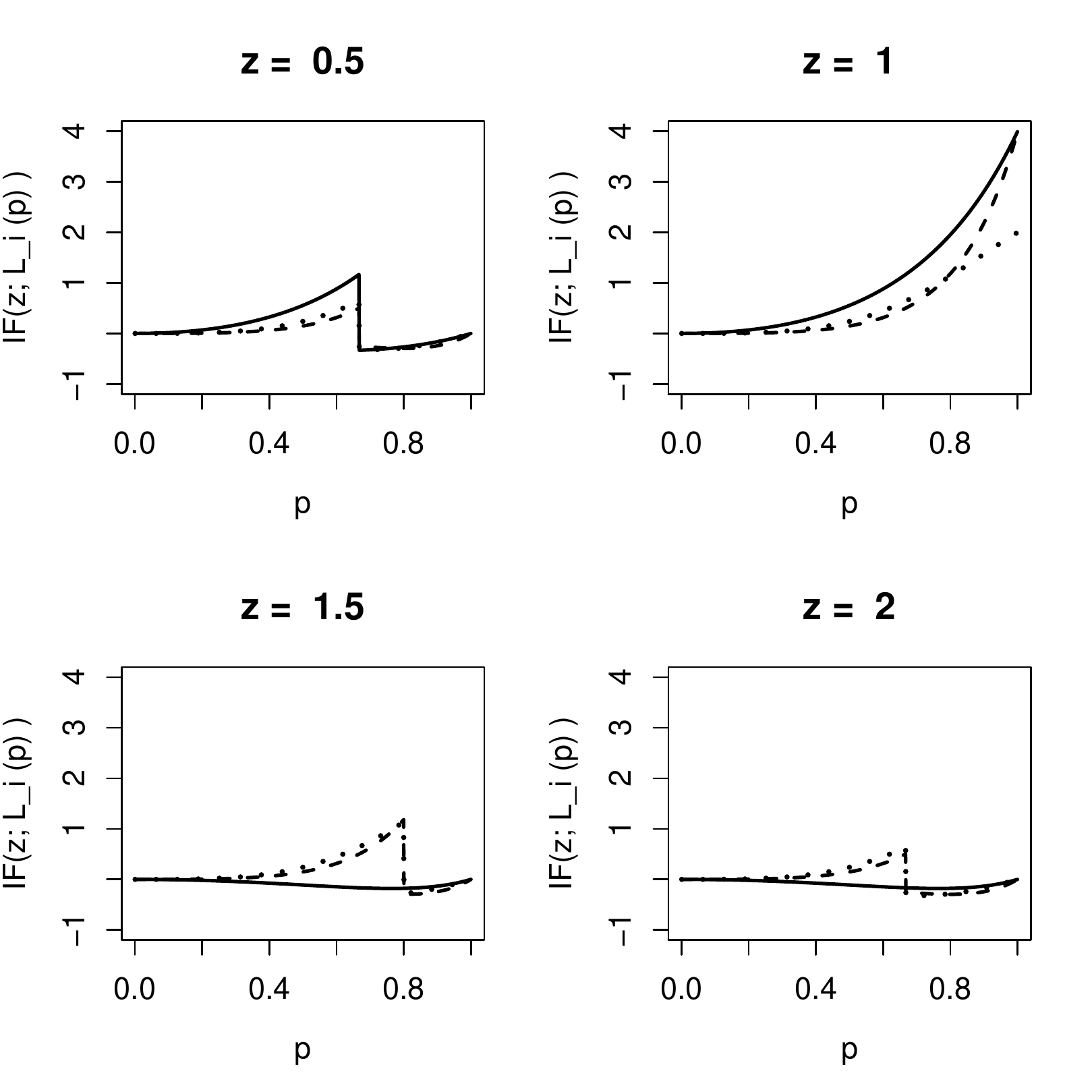}
\vspace{-.4cm}
\caption{\small \em  For various choices of $z$, $\IF (z;\,L _i(p),F_1)$  is plotted as a function of $p$.\label{fig7}}
\end{center}
\end{figure}

\clearpage
\newpage
\subsection{Influence functions of quantile inequality curves}\label{sec:IFqorenz}

\cite{cowell-2002} show that the influence function of the Lorenz curve at the point $p$ is unbounded, implying that a small amount of contamination can lead to a large bias in estimation of its value;
on the other hand the quantile inequality curves proposed here all have bounded influence functions, provided only that $F\in \F '.$
To see this, note that each $T_i(F)=L _i(F;p)=p x_{p/2}/d_i(p),$ where
$d_1(p)=x_{1/2}$, $d_2(p)=x_{1-p/2}$ and $d_3(p)=(x_{p/2}+x_{1-p/2})/2$ are
all quantile functionals or an average of them.

\begin{proposition}\label{P2}
The influence function of the functional defined by $T_i(F)=L _i(F;p)$ is a multiple $p$ of the derivative of the ratio of two functionals, so
by elementary calculus we have for each $p\in (0,1)$
\begin{equation}\label{IFLis} \nonumber
  \IF (z;\,T_i,F) = p\left \{\frac {\IF (z;\,x_{p/2},F)}{d_i(p)}
         -\frac {x_{p/2}\IF (z;\,d_i(p),F)}{d_i^2(p)} \right \}~.
 \end{equation}
For each case $i=1,2$ and 3 one only requires substitution of
the respective quantile influence functions for the $d_i$s found
in (\ref{IFquantile}).
\end{proposition}
 While these influence functions are complicated, the are easy to compute and plot
using currently available software. Specific examples are shown Figure~\ref{fig6} when the underlying $F=F_a$ is the Pareto distribution with shape parameter $a=1$ and are plotted as functions of a possible contamination at $z.$

  For small $p$ there is very little influence on $L _i(F;p)$ of contamination at any point $z$.  However, as $p$ increases, there is a noticeable increase in influence on $L _1(F_a;p)$ for $z$ near the median, which equals one in this case. Contamination at $z$ near zero is slightly negative, then rises to a positive relatively large positive peak  as $z$ approaches the median, and then drops to a small negative and constant influence again as $z$ increases past the median. This is to be  expected, because when the median is pulled to the left by contamination, then $L _1(F;p)=p\,x_{p/2}/x_{0.5}$ is increased, but when the median is pulled to the right,  the values of $L _1(F;p)$ are decreased.

 The other two $L _i(F;p)$ are similarly affected by contamination at $z$,
 but to a lesser extent.  Plots of the influence functions of the quantile inequality
 curves for other Pareto($a$) distributions (not shown) are similar to those in Figure~\ref{fig6}, and again the peak is located at the median  $F_a^{-1}(0.5)=2^{1/a}-1$.  Similar influence function plots are obtained  for uniform, lognormal and Weibull distributions, again with peaks
 near their respective medians.

\subsection{Influence of contamination at $z$ on the graph $\{p,L _i(p)\}$}

We have found, for each fixed $0<p<1$, the influence functions $\IF (z;\,L _i(p),F)$. Now we consider, for fixed $z$, the graph $\{(p,\IF (z;\,L _i(p),F))\}$, which shows the influence of contamination at $z$ on the respective inequality curves $\{(p,L _i(p))\}.$ Examples are shown in Figure~\ref{fig7}, again for $F$ the Pareto ($a=1$) distribution, and selected values of $z$.

First we concentrate on only the solid lines corresponding to $L _1(p)$. Inspection of (\ref{IFgini}) shows that the discontinuity points are $x_{1/2}=1$ and $x_{p/2}$.
Now $z<x_{p/2}$ if and only if $p>2F(z)$. Thus in the upper left plot of Figure~\ref{fig7}
where $z=0.5< x_{1/2}$ there are only two cases of interest: $p<2F(0.5)=2/3$ and $p>2/3$; in
the first interval $(0,2/3)$ the influence of contamination at $z=0.5$ on the $L _1$-curve
is positive and increasing in $p$, but its influence is negative for $p$ in $(2/3,1).$
For the top right plot $2F(z)=1$ so  the influence of contamination
$z=1$ at the median on the $L _1$-curve is positive and increasing for all $p$.
For the other two plots $z$ exceeds the median $2F(z)>1$ and there is only a slight
negative influence of $z$ on the $L _1$-curve for all $p$.

The influence of contamination at $z$ on the graphs of  $L _2(p)$, $L _3(p)$
is also shown in Figure~\ref{fig7} as dashed and dotted lines, respectively. Such influence is similar to that on $L _1(p)$ in the top two plots where $z$ does not exceed the median.
But in the lower plots where $z$ exceeds the median, the contamination is positive and
increasing on $(0,2(1-F(z)))$ and negative for larger $p$. For the bottom left plot this
interval is $(0,0.952)$, and for the bottom right it is $(0,0.8)$. Details are left as an
exercise for the reader.  Further increasing the values of $z$ only diminishes its effect on the graphs.

\begin{figure}[t!]
\begin{center}
\includegraphics[scale=.5]{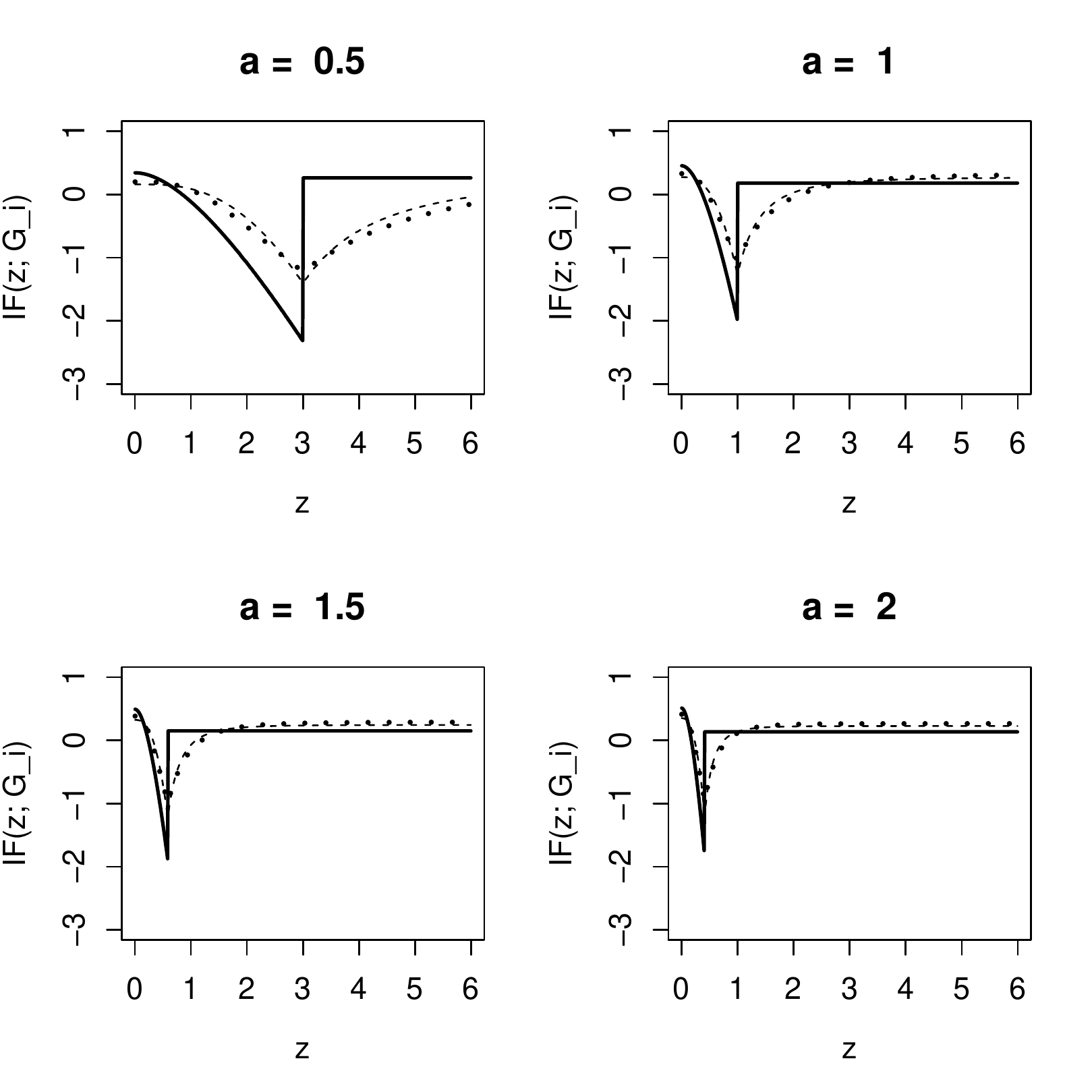}
\caption{\small \em The solid, dashed and dotted lines correspond, respectively, to the influence functions of $G _i$ for Pareto$(a)$ distributions for $i=1,2$ and 3.\label{fig8}}
\end{center}
\end{figure}

\subsection{Influence functions of quantile coefficients of inequality}\label{sec:IFgini}

The influence functions of the inequality coefficients associated with the $L _i$-curves are easily found, because the functional $G _i(F)=1- 2 \int _0^1L _i(F;p)\,dp $, which contains an average of $L _i(F;p)$ values over $p \in (0,1)$.

\begin{proposition}\label{P3}
For each $i=1,2$ and 3 the influence function of the inequality coefficients
$G_i$ are given respectively by
\begin{equation}\label{IFgini}
     \IF (z;G _i,F) =  -2 \int _0^1\IF (z;\,L _i(\,\cdot \,;p),F)\,dp ~.
\end{equation}
One only needs to justify taking the derivative $G _i(F_\epsilon ^{(z)})$
with respect to $\epsilon $ at $\epsilon =0$
under the integral sign. An argument based on the Leibniz Integration Rule
is given in the Appendix.
\end{proposition}
Figure~\ref{fig8} gives plots of the influence functions $\IF (z;\,G _i,F_a) =  -2 \int _0^1\IF (z;\,L _i(\,\cdot \,;p),F_a)dp $ of the inequality coefficients $G _i(F_a)$ when $F_a$ is the Pareto$(a)$ distribution for selected values of $a$. The biggest influence of contamination occurs at $z=F_a^{-1}(0.5)=2^{1/a}-1.$

The mean and variance of $\IF (z;G _i,F)$ are given by
\begin{eqnarray}\label{IFginimoments} \nonumber
 \e _F [\IF (Z;\,G _i,F)] &=&  -2 \int _0^1\e [\IF (Z;\,L _1(\;\cdot \,;p),F)]\;dp =0\\
 \var _F [\IF (Z;\,G _i,F)]&=&  4\;\e \left [\left \{\int _0^1\IF (Z;\,L _i(\,\cdot \,;p),F)\;dp \right \}^2\right ] ~.
 \end{eqnarray}
These quantities are easy to compute numerically;  examples
of the asymptotic standard error $\se [\hat G _i]= \{\var _F [\IF (Z;\,G _1,F)]\}^{1/2}$  determined by (\ref{IFginimoments}) are shown in Table~\ref{table2}.

\section{Convexity of  the quantile inequality curves}\label{sec:convex}

One of the nice mathematical properties of the Lorenz curve $\{p,L_0(F;p)\}$ is that it is convex for all distributions $F\in \F _0$.  The quantile-based versions (\ref{qcurves}) are defined for all $F $ in the larger class $\F $, but
need not be convex. In particular, empirical versions are often not convex over  $(0,1)$. The following examples
demonstrate that for the more commonly assumed income distributions, the quantile inequality curves are convex.
See \cite{J-K-B-1994,J-K-B-1995} for background material on these distributions.

\subsection{Non-convex example}\label{nonconvex}
\begin{figure}[t!]
\begin{center}
\includegraphics[scale=.6]{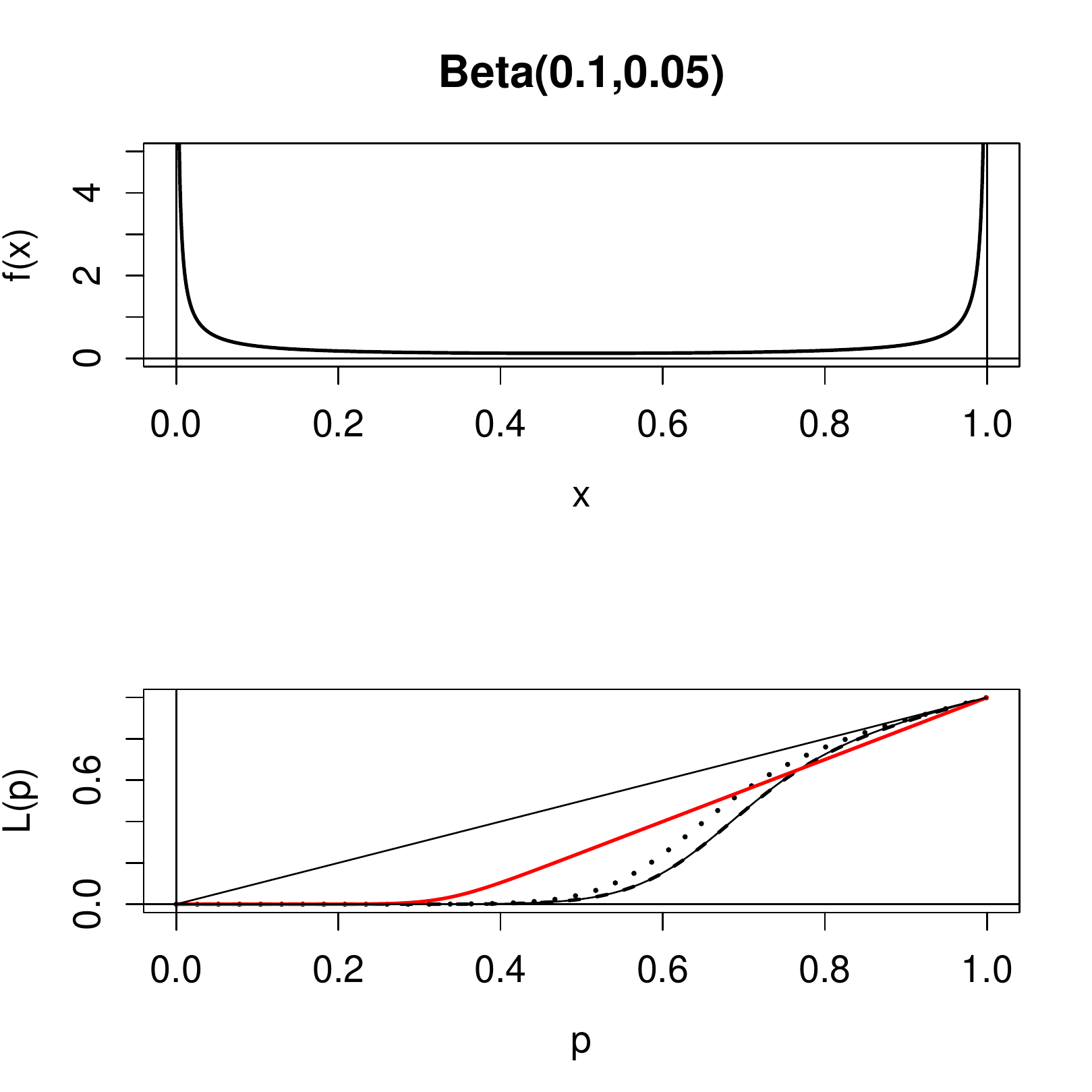}
\caption{\small \em The top plot shows the density of the Beta(0.1,0.05) distribution. Below it are the corresponding $L_i$ curves. The this solid, dashed and dotted lines correspond, respectively, to $i=1,2$ and 3;  The thick solid line is the Lorenz curve.\label{fig9}}
\end{center}
\end{figure}
Figure~\ref{fig9} shows that for the very U-shaped Beta distribution with parameters $(0.1, 0.05)$ only the Lorenz curve is convex. This distribution appears to have a symmetric density, but in fact is quite asymmetric, with mean  2/3, and the quartiles  0.050,0.997, and 1.000, to three decimal places.
The inequality coefficients are $G_0=0.329$, $G_1=0.453$, $G_2=0.455$ and $ G_3=0.403.$  Note that the Gini coefficient
$G_0< 1/3$, its value for the uniform distribution, a non-intuitive result to us.

Other plots, not shown, for parameters $(0.05,0.1)$, $(0.1,0.1)$ and $(0.05,0.05)$ indicate that
all four $L_i$ curves are convex.

\begin{table}[t!]
\begin{center}
\begin{small}
\begin{tabular}{lccccc}
                    &  $1-F(x)$            &   & $Q(p)$        & & $q(p)$                                 \\[.2cm]
\hline
 Exponential        &  $e^{-x} $           &   & $-\ln (1-p)$  & &  $(1-p)^{-1} $                         \\[.3cm]
 Normal             &  $\Phi (-x)$         &   & $z_p$         & & $\frac {1}{\varphi (z_p)}$             \\[.3cm]
 Lognormal          &  $ \Phi (-\ln (x))$  &   & $e^{z_p}$     & & $\frac {e^{z_p}}{\varphi(z_p)}$        \\[.3cm]
 Type I Pareto$(a)$ &  $ x^{-a} $          &   & $\frac{1}{(1-p)^{1/a}}$ & &$\frac{1}{a(1-p)^{1/a+1}}$     \\[.3cm]
Type II Pareto$(a)$ &  $ (1+x)^{-a} $      &   &$ \frac{1}{(1-p)^{1/a}}-1$ & &$\frac{1}{a(1-p)^{1/a+1}}$ \\[.3cm]
 Weibull$(\beta)$   &  $ e^{-x^{\beta}} $  &   &$\{-\ln (1-p)\}^{1/\beta}$ & &
 $\frac {\{-\ln (1-p)\}^{1/\beta -1}}{\beta(1-p)}$
\end{tabular}
\caption{{\bf Examples of distributions $F(x)$ and associated quantile functions and their densities.\ } \em In general, we denote $x_p=Q(p)=F^{-1}(p)$, but for the normal $F=\Phi $ with density $\varphi $, we write
$z_p=\Phi ^{-1}(p)$. the support of each $F$ is $(0,+\infty )$, except for the normal and Type I Pareto, the latter having
support on $[1,+\infty)$. \label{table4}}
\end{small}
\end{center}
\end{table}

\subsection{Convex examples}\label{convexexs}

\subsubsection*{Example 1.\quad Uniform.}

Starting with $Q(p)\equiv p$, we find $L_1(p)=p^2=L_3(p)$ and $L_2(p)=p^2/(2-p)$, all clearly convex functions of $p$ in (0,1).

\subsubsection*{Example 2.\quad Exponential.}

Here $Q(p)= -\ln(1-p)$, so $L_1(p)=-p\ln  \left( 1-p/2 \right)/\ln  \left( 2 \right)$ where $L_1''(p)=(4-p)/\left[ (p-2)^{2}\ln  \left( 2 \right)\right]>0$.  Similarly, $L_2(2)=p\ln  \left( 1-p/2 \right)/\ln  \left( p/2 \right)$ and $L_3(p)=2p\ln(1-p/2)/\ln[p(1-p/2)/2)$ and it is not difficult to show that both $L_2''(p)>0$ and $L_3''(p)>0$ so that $L_1(p)$, $L_2(p)$ and $L_3(p)$ are all convex.

\subsubsection*{Example 3.\quad Lognormal.}
It is \lq obvious\rq\  from the lower left plot in Figure~\ref{fig1} that all three $L _i(p)$ curves are convex on (0,1) for the lognormal distribution.  Proving it using the calculus is not as straightforward as one might expect. Note that $Q(p)= e^{z_p}$, $q(p)=e^{z_p}/\varphi (z_p).$ Further, observe that $L_1(p) =p\exp (z_{p/2})$ and that $\exp (z_{p/2})$ is not convex, so one cannot use the fact that two monotone increasing convex functions is convex.  Taking derivatives,
\begin{eqnarray*}
L_1'(p) &=&  L_1(p)\left \{ \frac {1}{p} +\frac{1}{2\varphi (z_{p/2})} \right \}  \\
L_1''(p) &=& L_1(p)\left [ \left \{ \frac {1}{p} +\frac{1}{2\varphi (z_{p/2})} \right \}^2-\frac{1}{p^2}-\frac{\varphi '(z_{p/2})}{4\varphi ^3 (z_{p/2})}         \right ] \\
       &=& L_1(p)\left [ \frac{1}{p\varphi (z_{p/2})} +\frac{1+z_{p/2}}{4\varphi ^2(z_{p/2})}\right ]  ~.
\end{eqnarray*}
Thus $L''_1(p)>0$ if and only if $4\varphi (z_{p/2})+p(1+z_{p/2})>0$ and this again, while obvious from a plot,
is not readily verified.

Next consider $L_2(p)=p\,\{\exp (z_{p/2})\exp (-z_{1-p/2})\}=p\,\exp (2z_{p/2})$. The argument is
very similar to that for $L_1$:
\begin{eqnarray*}
L_2'(p) &=&  L_2(p)\left \{ \frac {1}{p} +\frac{1}{\varphi (z_{p/2})} \right \}  \\
L_2''(p) &=& L_2(p)\left [ \left \{ \frac {1}{p} +\frac{1}{\varphi (z_{p/2})} \right \}^2-
\frac{1}{p^2}-\frac{\varphi '(z_{p/2})}{2\varphi ^3 (z_{p/2})} \right ] \\
       &=& L_2(p)\left [ \frac{2}{p\varphi (z_{p/2})} +\frac{2+z_{p/2}}{2\varphi ^2(z_{p/2})}\right ]  ~.
\end{eqnarray*}
Thus $L''_2(p)>0$ if and only if $4\varphi (z_{p/2})+p(2+z_{p/2})>0$, a weaker condition than required for
convexity of $L_1$.

Finally, consider $L_3(p)=2p/\{1+p/L _2(p)\}=2p/\{1+\exp (-2z_{p/2})\}$. It suffices to show that
$h(p)=1/\{1+\exp (-2z_{p/2})\}$ is convex in $p$ and this is readily verified.

\subsubsection*{Example 4.\quad Type I Pareto.}

For the Type I Pareto$(a)$ distribution where $a>0$, $Q(p)=(1-p)^{-1/a}$.  Let $c_1=(2-p)^{-1/a}/a$ which is positive.  Then $L_1''(p)=c_1[(1-p/2)^{-1}+(1+1/a)p/(p-2)^2]>0$ so that $L_1(p)$ is convex.  Similarly, $L_2''(p)=c_1p^{1/a}(1+1/a)[(1-p/2)^{-1}+p/(p-2)^2+1/p]>0$ so that $L_2(p)$ is also convex.  The expression for $L_3''(p)$ is much more complicated although plots and computational minimization reveal that convexity holds.  For example, over all $p\in [0,1)$ and $a\in (0, 10]$, min $L_3''(p)=0.169$ (at $p=0.667$ and $a=10$).

\subsubsection*{Example 5.\quad Type II Pareto.}

For the Type II Pareto$(a)$ distribution where $a>0$, $Q(p)=(1-p)^{-1/a}-1$. We then have that
\begin{align*}
L_1''(p)=&\frac{(1-p/2)^{-1/2}}{a^2(p-2)^2(2^{1/a}-1)}\left[p + a(4-p)\right]>0
\end{align*}
so that $L_1(p)$ is convex.  Both $L_2''(p)$ and $L_3''(p)$ are complicated expressions although computational minimization reveals non-negative minimums over all $p$ and $a\in(0,10]$.

\subsubsection*{Example 6.\quad Weibull.}

For the Weibull distribution with shape parameter $\beta>0$, we have
$$L_1''(p)=\frac{\ln(2)^{-1/\beta}}{\beta(p-2)^2}\ln\left(\frac{2}{2-p}\right)^{1/\beta-1}\left[4-p-p\ln\left(\frac{2}{2-p}\right)^{-1}+\frac{p}{\beta}\ln\left(\frac{2}{2-p}\right)^{-1}\right].$$
The term $-p\ln(2/(2-p))$ is a decreasing function in $p$ with limit equal to -2 as $p$ approaches 0.  Consequently, $L_1''(p)>0$ so that $L_1(p)$ is convex.  For $L_2(p)$ and $L_3(p)$, again we used computational minimization for all $\beta$ values up to 100.  Neither had a negative minimum so both were found to be convex.

\section{Summary and further research}\label{sec:summary}

We have shown that quantile versions of the Lorenz curve have most of  the properties of the original definition, with two exceptions. The first exception
is convexity, which is not satisfied for some very U-shaped distributions and
many empirical ones. Nevertheless,
for all continuous distributions commonly used to model population incomes, the
quantile versions are convex.

The second exception is the first order transference principle, which is mean-centric. When replaced by a median-centric definition, this principle is
satisfied for all three quantile versions of the Lorenz curve.
We then studied a specific example of a transfer function and showed how it
could be measured by the associated inequality coefficients, defined as
twice the area between the quantile inequality curve and the equity diagonal.
These inequality coefficients can also be interpreted as expected values of
certain functions of independent randomly drawn incomes from the population.

These concepts have distinct advantages over the traditional Lorenz curve and
Gini index.  They are defined for {\em all} positive income distributions, and
their influence functions are bounded. Distribution-free confidence intervals for the ordinates of inequality curves at fixed points are readily found, since they are just ratios of finite linear combinations
of quantiles. In addition, we showed that the standard errors of estimates for the quantile analogues of the Gini coefficient do not appear to depend much on the underlying population model, so that sample sizes can be chosen in advance to obtain desired standard
errors. Simulation studies suggest that these sample inequality coefficients approach normality very rapidly, and confidence intervals for them can be constructed when the
underlying scale family is known. One way to
obtain distribution-free confidence intervals for them would be to find distribution-free estimates of their standard errors, which involves quantile density estimation.

Many other challenges remain.  It would be good to have simple necessary
and sufficient conditions in terms of the underlying income distribution
for convexity of each of the inequality curves.
 If one is interested in confidence bands for the quantile curves, one could
utilize functionals of the quantile process to determine them, starting with
the results in \cite{doss-1992}.
Finally, applications to other fields which use diversity indices \cite{patil-1982} would be of interest, as well as
connections to the \lq Lorenz dominance\rq\  literature, see \cite{aab-2011} and
references therein.

\section{Appendix:\quad Proof of Proposition~3}

The interchange of limit (as $\epsilon \downarrow 0$) and integral is justified by the Leibniz Integral Rule. It requires that $h_i(p)\equiv \IF (z;\,L _i(\,\cdot \,;p),F)$ be continuous in $p$,
and bounded in absolute value for  $p\in (0,1)$ by an integrable function.

\subsubsection*{Proof for $i=1$.}

For $L _1$, we have from Proposition~\ref{P2} that
\begin{eqnarray}\label{IFL1bound}\nonumber
  |h_1(p)|&\leq & \frac {p}{x_{1/2}^2}\left \{x_{1/2}| \IF (z;\,Q(\,\cdot \, ,p/2),F)|
         +x_{p/2}|\IF (z;\,Q(\,\cdot\, ;1/2),F)| \right \}\\ \nonumber
 &\leq &  \frac {p}{x_{1/2}^2}\left \{x_{1/2}\;\max \{p/2,1-p/2\} q(p/2)
         +\frac {x_{p/2}\,q(1/2)}{2} \right \}~.
\end{eqnarray}
The second term is bounded because $p\,Q(p/2)\leq x_{1/2}$ for $p\in (0,1)$\; ; and, for the first term we require only that $p\,q(p/2)$  be integrable on $(0,1)$.
By making the change of variable $x=F^{-1}(p/2)$ in $\int _0^1p\,q(p/2)\,dp $ one
finds that this integral is bounded  by $4x_{1/2}.$ Therefore $|h_1(p)|$ is bounded
by an integrable function on $(0,1)$, justifying (\ref{IFgini}) for $L _1$.

\subsubsection*{Proof for $i=2$.}

 For $L _2(p)=p\,x_{p/2}/x_{1-p/2}$ we have
\begin{eqnarray}\label{IFL2bound}\nonumber
  h_2(p)&\equiv & \frac {p}{x_{1-p/2}^2}\,\left \{x_{1-p/2}\; \IF (z;\,Q(\,\cdot \, ,p/2),F)-
         x_{p/2}\;\IF (z;\,Q(\,\cdot\, ;1-p/2),F)\right \} \text { , so} \\
  |h_2(p)|&\leq &  \frac {p\,q(p/2)}{x_{1-p/2}}+\frac {p\,x_{p/2}\,q(1-p/2)}{x_{1-p/2}^2} ~.
\end{eqnarray}
 The first term in the last line of (\ref{IFL2bound}) is bounded above by $p\,q(p/2)/x_{1/2},$
and it has already been shown that $p\,q(p/2)$ was integrable on (0,1).

Next we show that the second term is bounded by an integrable function. Let $m=x_{1/2}$
and make the change of variable $x=F^{-1}(1-p/2)=x_{1-p/2}$ to obtain:
\begin{eqnarray}\label{IFL2bound2}
\nonumber
\int _0^1\frac {p\,x_{p/2}\,q(1-p/2)}{x_{1-p/2}^2}\,dp  &=&4\int _m^\infty
\frac {\{1-F(x)\}\;F^{-1}(1-F(x))}{x^2}\,dx \\
   &\leq &  4m\int _m^\infty  \frac{dx}{x^2} =4~.
\end{eqnarray}
This shows that $h_2(p)= \IF (z;\,L _2(\,\cdot \,;p),F)$ is bounded on  $(0,1)$ by an integrable function.
\subsubsection*{Proof for $i=3$.}
Let $m(p)=(x_{p/2}+x_{1-p/2})/2,$ so $m(1)=m$ is the median, and $L _3(p)=p\,x_{p/2}/m(p).$
It is immediate that $\IF (z;\,m(p),F)=\{\IF (z;\,Q(\,\cdot \, ,p/2),F)+\IF (z;\,Q(\,\cdot \, ,1-p/2),F)\}/2$ and that $|\IF (z;\,m(p),F)|\leq \{q(p/2)+q(1-p/2)\}/2.$

Consider bounding $h_3(p)= \IF (z;\,L _3(\,\cdot \,;p),F)$ by an integrable function.
\begin{eqnarray}\label{IFL3bound}\nonumber
  h_3(p)&\equiv & \frac {p}{m^2(p)}\,\left \{m(p)\; \IF (z;\,Q(\,\cdot \, ,p/2),F)-
         x_{p/2}\;\IF (z;\,m(p),F)\right \} \text { , so} \\
  |h_3(p)|&\leq &  \frac {p\,q(p/2)}{m(p)}+\frac {p\,x_{p/2}\,\{q(p/2)+q(1-p/2)\}}{2m^2(p)} ~.
\end{eqnarray}
The first term $p\,q(p/2)/m(p)\leq 2p\,q(p/2)/x_{1-p/2}$, which has already shown to be integrable.
The third term $p\,x_{p/2}\,q(1-p/2)/(2m^2(p))\leq 2p\,x_{p/2}\,q(1-p/2)/x_{1-p/2}^2,$ shown
to be integrable in (\ref{IFL2bound2}).  The second term $p\,x_{p/2}\,q(p/2)/(2m^2(p))\leq p\,q(p/2)/x_{1-p/2},$ using the fact that $m^2(p)\geq x_{p/2}x_{1-p/2}$.
Therefore $|h_3(p)|$ is bounded by an integrable function.


\begin{thebibliography}{}

\bibitem[\protect\citename{Aaberge \& Mogstad, }2011]{aab-2011}
{\sc Aaberge, M., \& Mogstad, R.} 2011.
\newblock Robust inequality measures.
\newblock {\em Journal of {E}conomic {I}nequality}, {\bf 9}(3), 353--371.

\bibitem[\protect\citename{Beach \& Davidson, }1983]{beach-1983}
{\sc Beach, C.M., \& Davidson, R.} 1983.
\newblock Distribution-free statistical inference with lorenz curves and income
  shares.
\newblock {\em Review of {E}conomic {S}tudies}, {\bf L}, 723--735.

\bibitem[\protect\citename{Brown, }1981]{brown-1981}
{\sc Brown, B.M.} 1981.
\newblock Symmetric quantile averages and related estimators.
\newblock {\em Biometrika}, {\bf 68}(1), 235--242.

\bibitem[\protect\citename{Cowell \& Victoria-Feser, }1996]{cowell-1996}
{\sc Cowell, F.A., \& Victoria-Feser, M.P.} 1996.
\newblock Robustness properties of inequality measures.
\newblock {\em Econometrica}, {\bf 64}(1), 77--101.

\bibitem[\protect\citename{Cowell \& Victoria-Feser, }2002]{cowell-2002}
{\sc Cowell, F.A., \& Victoria-Feser, M.P.} 2002.
\newblock Welfare rankings in the presence of contaminated data.
\newblock {\em Econometrica}, {\bf 70}(3), 1221--1233.

\bibitem[\protect\citename{Cowell \& Victoria-Feser, }2003]{cowell-2003}
{\sc Cowell, F.A., \& Victoria-Feser, M.P.} 2003.
\newblock Distribution-free inference for welfare indices under complete and
  incomplete information.
\newblock {\em The {J}ournal of {E}conomic {I}nequality}, {\bf 1}(3), 191--219.

\bibitem[\protect\citename{Cowell \& Victoria-Feser, }2007]{cowell-2007}
{\sc Cowell, F.A., \& Victoria-Feser, M.P.} 2007.
\newblock Robust stochastic dominance: A semi-parametric approach.
\newblock {\em The {J}ournal of {E}conomic {I}nequality}, {\bf 5}(1), 21--37.

\bibitem[\protect\citename{Dalton, }1920]{dalton-1920}
{\sc Dalton, H.} 1920.
\newblock The measurement of the inequality of incomes.
\newblock {\em Economic {J}ournal}, {\bf 30}, 348--361.

\bibitem[\protect\citename{Davidson, }2008]{davidson-2008}
{\sc Davidson, R.} 2008.
\newblock Reliable inference for the {G}ini index.
\newblock {\em Journal of {E}conometrics}, {\bf 150}, 30--40.

\bibitem[\protect\citename{{D}evelopment~{C}ore {T}eam, }2008]{R}
{\sc {D}evelopment~{C}ore {T}eam, R}. 2008.
\newblock {\em R: A language and environment for statistical computing}.
\newblock R {F}oundation for {S}tatistical {C}omputing, Vienna, Austria.
\newblock {ISBN} 3-900051-07-0.

\bibitem[\protect\citename{Doss \& Gill, }1992]{doss-1992}
{\sc Doss, H., \& Gill, R.D.} 1992.
\newblock An elementary approach to weak convergence for quantile processes.
\newblock {\em Journal of the {A}merican {S}tatistical {A}ssociation}, {\bf
  87}, 869--877.

\bibitem[\protect\citename{Fellman, }2012]{fell-2012}
{\sc Fellman, J.} 2012.
\newblock Properties of {L}orenz curves for transformed income distributions.
\newblock {\em Theoretical {E}conomics {L}etters}, {\bf 2}, 487--493.

\bibitem[\protect\citename{Gastwirth, }1971]{gast-1971}
{\sc Gastwirth, J.L.} 1971.
\newblock A general definition of the {L}orenz curve.
\newblock {\em Econometrika}, {\bf 39}, 1037--1039.

\bibitem[\protect\citename{Gastwirth, }2012]{gast-2012}
{\sc Gastwirth, J.L.} 2012.
\newblock A robust {G}ini-type index better detects in the income distribution:
  a reanalysis of income distribution in the {U}nited {S}tates from 1967-2011.
\newblock {\em {SSRN} {E}lectronic {J}ournal}.
\newblock {DOI}: 10.2139/ssrn.2164745.

\bibitem[\protect\citename{Gini, }1914]{gini-1914}
{\sc Gini, C.} 1914.
\newblock Sulla misura della concentrazione e della variabilit`a dei caratteri.
\newblock {\em {A}tti del {R}eale {I}stituto {V}eneto di {S}cienze, {L}ettere
  ed {A}rti}, {\bf 73}, 1203--1248.
\newblock English translation (2005) in {\em {M}etron} Vol. 63, pp. 3--38.

\bibitem[\protect\citename{Hampel {\em et~al.}, }1986]{HRRS86}
{\sc Hampel, F.R., Ronchetti, E.M., Rousseeuw, P.J., \& Stahel, W.A.} 1986.
\newblock {\em Robust {S}tatistics: The {A}pproach {B}ased on {I}nfluence
  {F}unctions}.
\newblock New York: John Wiley and Sons.

\bibitem[\protect\citename{Hyndman \& Fan, }1996]{hynd-1996}
{\sc Hyndman, R.J., \& Fan, Y.} 1996.
\newblock Sample quantiles in statistical packages.
\newblock {\em The {A}merican {S}tatistician}, {\bf 50}, 361--365.

\bibitem[\protect\citename{Johnson {\em et~al.}, }1994]{J-K-B-1994}
{\sc Johnson, N.L., Kotz, S., \& Balakrishnan, N.} 1994.
\newblock {\em Continuous {U}nivariate {D}istributions}.
\newblock  Vol. 1.
\newblock New York: John Wiley \& Sons.

\bibitem[\protect\citename{Johnson {\em et~al.}, }1995]{J-K-B-1995}
{\sc Johnson, N.L., Kotz, S., \& Balakrishnan, N.} 1995.
\newblock {\em Continuous {U}nivariate {D}istributions}.
\newblock  Vol. 2.
\newblock New York: John Wiley \& Sons.

\bibitem[\protect\citename{Kampke, }2010]{kampke-2010}
{\sc Kampke, T.} 2010.
\newblock The use of mean values vs. medians in inequality analysis.
\newblock {\em Journal of {E}conomic and {S}ocial {M}easurement}, {\bf 35},
  43--62.

\bibitem[\protect\citename{Kleiber, }2005]{kleiber-2005}
{\sc Kleiber, C.} 2005.
\newblock {\em The {L}orenz curve in {E}conomics and {E}conometrics}.
\newblock Technical Report~30. University of {D}ortmund, SFB 475.

\bibitem[\protect\citename{Parzen, }1979]{par-1979}
{\sc Parzen, E.} 1979.
\newblock Nonparametric statistical data modeling.
\newblock {\em Journal of the {A}merican {S}tatistical {A}ssociation}, {\bf 7},
  105--131.

\bibitem[\protect\citename{Patil \& Taillie, }1982]{patil-1982}
{\sc Patil, G.P., \& Taillie, C.} 1982.
\newblock Diversity as a concept and its measurement.
\newblock {\em Journal of the {A}merican {S}tatistical {A}ssociation}, {\bf
  77}, 548--561.

\bibitem[\protect\citename{Prendergast \& Staudte, }2015a]{ps-2015a}
{\sc Prendergast, L.A., \& Staudte, R.G.} 2015a.
\newblock {Exploiting the Quantile Optimality Ratio to Obtain Better Confidence
  Intervals of Quantiles}.
\newblock {\em {arXiv preprint arXiv:1505.04234}}.

\bibitem[\protect\citename{Prendergast \& Staudte, }2015b]{ps-2015b}
{\sc Prendergast, L.A., \& Staudte, R.G.} 2015b.
\newblock {When large n is not enough-Distribution-free Interval Estimators for
  Ratios of Quantiles}.
\newblock {\em ar{X}iv preprint ar{X}iv:1508.06321v2}.

\bibitem[\protect\citename{Sen, }1986]{sen-1986}
{\sc Sen, P.K.} 1986.
\newblock The {G}ini coefficient and poverty indexes: some reconciliations.
\newblock {\em Journal of the {A}merican {S}tatistical {A}ssociation}, {\bf
  81}, 1050--1057.

\bibitem[\protect\citename{Sheather \& Marron, }1990]{shma-1990}
{\sc Sheather, S.J., \& Marron, J.S.} 1990.
\newblock Kernel quantile estimators.
\newblock {\em Journal of the {A}merican {S}tatistical {A}ssociation}, {\bf
  85}, 410--416.

\bibitem[\protect\citename{Staudte, }2013b]{S-2013}
{\sc Staudte, R.G.} 2013b.
\newblock Distribution-free confidence intervals for the standardized median.
\newblock {\em {STAT}}, {\bf 2}(1), 184--196.

\bibitem[\protect\citename{Staudte, }2014]{S-2014}
{\sc Staudte, R.G.} 2014.
\newblock Inference for quantile measures of skewness.
\newblock {\em {Test}}, {\bf 23}(4), 751--768.

\bibitem[\protect\citename{Staudte, }2015]{S-2015}
{\sc Staudte, R.G.} 2015.
\newblock Inference for quantile measures of kurtosis, peakedness and
  tail-weight.
\newblock {\em Communications in {S}tatistics}.
\newblock In press.

\bibitem[\protect\citename{Staudte \& Sheather, }1990]{S-S-1990}
{\sc Staudte, R.G., \& Sheather, Simon~J.} 1990.
\newblock {\em Robust {E}stimation and {T}esting}.
\newblock New York: Wiley.

\bibitem[\protect\citename{Tukey, }1965]{tukey-1965}
{\sc Tukey, J.W.} 1965.
\newblock Which part of the sample contains the information?
\newblock {\em Proceedings of the {M}athemetical {A}cademy of {S}cience {USA}},
  {\bf 53}, 127--134.

\bibitem[\protect\citename{Victoria-Feser, }2000]{VF-2000}
{\sc Victoria-Feser, M.P.} 2000.
\newblock Robust methods for the analysis of income distribution, inequality
  and poverty.
\newblock {\em International {S}tatistical {R}eview}, {\bf 68}(3), 277--293.

\bibitem[\protect\citename{Victoria-Feser \& Ronchetti, }1994]{VFronc-1994}
{\sc Victoria-Feser, M.P., \& Ronchetti, E.} 1994.
\newblock Robust methods for personal-income distribution models.
\newblock {\em The {C}anadian {J}ournal of {S}tatistics}, {\bf 22}(2),
  247--258.

\end{thebibliography}
\end{document}